\documentstyle[11pt]{article}

\newcommand{\C}{ {\bf C} }
\newcommand{\Z}{ {\bf Z} }
\newcommand{\N}{ {\bf N} }
\newcommand{\Scal}{ {\cal S} }
\newcommand{\M}{ {\cal M} }
\newcommand{\G}{ {\cal G} }
\newcommand{\A}{ {\cal A} }
\newcommand{\f}{ \widetilde{F} }
\newcommand{\sgp}{ {\cal S}_{G,P} }
\newcommand{\Gzero}{ {\cal G}^{(o)} }
\newcommand{\Szero}{ {\cal S}^{(o)} }
\newcommand{\lciup}{ \widehat{L} }
\newcommand{\dnciup}{ \widehat{D_{n}} }
\newcommand{\tilbeta}{ \widetilde{\beta} }
\newcommand{\pq}{ PP^{-1} }
\newcommand{\xprod}{ x_{1} \cdots x_{n} }
\newcommand{\xuri}{ x_{1}, \ldots , x_{n} }
\newcommand{\igreci}{ y_{1} \cdots  y_{n} }
\newcommand{\beturi}{ \beta_{x_{1}} \cdots \beta_{x_{n}} }
\newcommand{\ecdef}{ \stackrel{def}{\Leftrightarrow} }
\newcommand{\invsg}{ \mbox{inverse semigroup} }
\newcommand{\siminvsg}{ \mbox{symmetric inverse semigroup} }
\newcommand{\local}{ \mbox{localization} }
\newcommand{\zem}{ \mbox{zero element} }
\newcommand{\homo}{ \mbox{homomorphism} }
\newcommand{\homeo}{ \mbox{homeomorphism} }

\begin{document}

\title{\bf On a groupoid construction for actions of certain
inverse semigroups}
\author{Alexandru Nica  \\
Department of Mathematics  \\
University of California  \\
Berkeley, California 94720    \\
(e-mail: nica@cory.berkeley.edu)}
\date{ }

\maketitle

\vspace{1.5in}

\begin{abstract}
We consider a version of the notion of $F$-inverse semigroup (studied
\setcounter{page}{0}
in the algebraic theory of inverse semigroups). We point out that an
action of such an inverse semigroup on a locally compact space has
associated a natural groupoid construction, very similar to the
one of a transformation group. We discuss examples related to
Toeplitz algebras on subsemigroups of discrete groups, to
Cuntz-Krieger algebras, and to crossed-products by partial
automorphisms in the sense of Exel.
\end{abstract}

\vspace{3in}
\setlength{\baselineskip}{18pt}

\setcounter{section}{-1}
\section{Introduction} \

Let $G$ be a discrete group, let $P$ be a unital subsemigroup
of $G$, and let ${\cal W} (G,P)$ denote the Toeplitz (also
called Wiener-Hopf) C*-algebra associated to $(G,P)$; i.e.,
${\cal W}(G,P) \subseteq {\cal L}(l^{2}(P))$
is the C*-algebra generated by the compression to $l^{2}(P)$ of
the left regular representation $\Lambda : l^{1}(G) \rightarrow
{\cal L}(l^{2}(G))$.

A powerful tool for studying ${\cal W}(G,P)$ is a locally
compact groupoid $\G$ introduced by P. Muhly and J. Renault
in \cite{MR} (in a framework larger than the one of discrete
groups), and which was shown in \cite{MR} to have
${\mbox{C*}}_{red} ( \G ) \simeq {\cal W}(G,P)$. We shall
refer to $\G$ as to the {\em Wiener-Hopf groupoid} associated
to $(G,P)$.

On the other hand let us denote, for every $x \in G$, by
$\beta_{x} : \{ t \in P \ | \ xt \in P \}$
$\rightarrow$ $\{ s \in P \ |$ $x^{-1}s \in P \}$ the
partially defined left translation with $x$ on $P$; and
let us denote by $\sgp$ the semigroup of bijective
transformations between subsets of $P$ (with multiplication
given by composition, defined where it makes sense) which is
generated by $( \beta_{x} )_{x \in G}$. Then $\sgp$ is an
{\em $\invsg$}, which encodes in some sense the action of $G$
on $P$ by left translations.

The original motivation for this work was to understand
the relation between the $\invsg$ $\sgp$ and the Wiener-Hopf
groupoid associated to $(G,P)$. It turns out that:

$\ \ \ \ $(A) there exists a remarkable action of $\sgp$ on a compact
space,

$\ \ \ \ $and

$\ \ \ \ $(B) there exists a very natural procedure of associating a
locally compact groupoid to an
action of an $\invsg$ in a class which contains $\sgp$,

$\ \ \ $
\linebreak
such that (A) and (B) together lead to the Wiener-Hopf
groupoid.

Moreover, there exists a natural unital $\star$-homomorphism
$\Psi$ from C*($\sgp$) onto C*($\G$), where C*($\sgp$) denotes
the enveloping C*-algebra of $\sgp$, and C*($\G$) the
C*-algebra of the Wiener-Hopf groupoid. Verifying whether $\Psi$
is faithful can be reduced to studying the surjectivity of a
map between compact spaces, and actually to comparing two
subspaces of ${\{ 0,1 \}}^{G}$ (see Corollary 6.4 and
Example 6.5 below). This holds when $G$, viewed as a
left-ordered group with positive semigroup $P$, has some sort
of lattice properties (as for instance those considered in
\cite{N2}), but is not true in general.

Let us now return to the above (A) and (B). By an action of
an $\invsg$ $\Scal$ on a locally compact space $\Omega$ we
shall understand a unital $\star$-homomorphism from $\Scal$
into the $\invsg$  of $\homeo$s between open subsets of
$\Omega$.

Related to (A), let us assume for a moment that
$P$ does not contain a proper subgroup of $G$; then the
action of $\sgp$ mentioned at (A), let us call it $\Phi$,
is on a space $\Omega$ which is a compactification of $P$.
Every $\alpha \in \sgp$ can thus be viewed as a
function between two subsets $Dom( \alpha )$,
$Ran( \alpha ) \subseteq P \subseteq \Omega$, and we have the
remarkable fact that $\Phi ( \alpha ) $ is the unique
continuous extension of $\alpha$ between the closures
$\overline{Dom( \alpha )}$ and
$\overline{Ran( \alpha )}$ in $\Omega$.

Related to (B): the fact which makes the groupoid
construction associated to an action of $\sgp$ to be very
simple, and indeed a straightforward generalization of
the groupoid associated to a transformation group, is
the following: every element of $\sgp$, except possibly
the zero element, is majorized ( in the sense of the usual
partial order on an $\invsg$) by a unique maximal element.
Inverse semigroups with this property have been studied
in the algebraic theory of $\invsg$s under the name of
$F$-$\invsg$s (see \cite{Pet}, Section VII.6).

There are various examples of $F$-$\invsg$s coming from the
algebraic theory which can be considered (see Example 1.4
below). Given the nature of the present note, it is probably
even more interesting to remark that most of the singly
generated $\invsg$s have the property under consideration;
moreover, the C*-algebra of the groupoid associated to an
action of a singly generated $\invsg$ is isomorphic to the
corresponding crossed-product C*-algebra by a partial
automorphism, in the sense of R. Exel \cite{E1}.

In \cite{K}, A. Kumjian has developed a method which
associates a C*-algebra to a {\em localization}, i.e. to
an $\invsg$ $\Scal$ of $\homeo$s between
open subsets of a locally compact
space, such that the family of domains
$(Dom ( \alpha ))_{\alpha \in \Scal}$ is a basis for the
topology of the space. The point of view of the present note
is slightly different from the one of \cite{K}, and the
examples we consider do not generally have the named
localization property (see Section 5 below). For $\Scal$ in
the common territory of the two approaches, the C*-algebra of
the associated groupoid can be shown, however, to be isomorphic
to the one constructed in \cite{K} (see Theorem 5.1 below).
This is for instance the case for the localizations
giving the Cuntz-Krieger C*-algebras ${\cal O}_{A}$.

The paper is subdivided into sections as follows: after
recalling in Section 1 some basic definitions, and giving
some examples, the groupoid construction is presented in
Section 2. In Section 3 we discuss the example related to
Toeplitz algebras, and in Section 4 the one related to
crossed products by partial automorphisms. Section 5 is
devoted to the relation with localizations; in 5.4 the example
related to the Cuntz-Krieger algebras is discussed. Finally,
in Section 6 we study the relation between C*($\Scal$) and
the C*-algebra of the associated groupoid (and obtain in
particular the facts about $\sgp$ mentioned
at the bottom of page 1).

$\ $

$\ $

\section{Basic definitions and examples} \

A semigroup $\Scal$ is an {\em $\invsg$} if for
every $\alpha \in \Scal$
there exists a unique element of $\Scal$, denoted by
$\alpha^{*}$, such that
$\alpha \alpha^{*} \alpha = \alpha$,
$\alpha^{*} \alpha \alpha^{*} = \alpha^{*}$. For the
algebraic facts about $\invsg$s needed in this note, we shall
use as a reference the monograph \cite{Pet}. For $T$ a
non-void set, we denote by ${\cal I}_{T}$ the
{\em symmetric inverse semigroup} on $T$, i.e. the
semigroup of all bijections between subsets of $T$ (with
multiplication given by composition, defined where it
makes sense). This is a generic example, in the sense that
every $\invsg$ can be embedded into an ${\cal I}_{T}$.
(see \cite{Pet}, IV.1). On an arbitrary $\invsg$ we have
a natural partial order, defined by
\begin{equation}
\alpha \leq \beta \  \ecdef
\beta^{*} \alpha = \alpha^{*}  \alpha , \ \
\alpha , \beta \in \Scal .
\end{equation}
$\alpha \leq \beta$ is equivalent to the fact that $\alpha$
is a restriction of $\beta$ in every embedding of $\Scal$ into
a $\siminvsg$ (see \cite{Pet}, II.1.6 and IV.1.10). An element
$\alpha \in \Scal$ will be called maximal if there is no
$\beta \neq \alpha$ in $\Scal$ such that $\alpha \leq \beta$.

We shall consider only unital $\invsg$s, and the unit will be
usually denoted by $\epsilon$. Also, an $\invsg$ may or may not
have a $\zem$, which, if existing, is unique and will be denoted
by $\theta$ ( $\alpha \theta = \theta \alpha = \theta$ for all
$\alpha \in \Scal$).
If $\alpha$ is an element of the $\invsg$ ${\cal S}$, the
expression ``$\alpha $ is not zero element for ${\cal S}$''
will be used to mean that, if ${\cal S}$ happens to have
a zero element $\theta$, then $\alpha \neq \theta$.

A unital $\invsg$ ${\cal S}$ is called an
{\em $F$-inverse semigroup} if every $\alpha \in {\cal S}$
is majorized (in the sense of (1.1)) by a unique maximal
element of ${\cal S}$ (see \cite{Pet}, Definition VII.5.13).
This is essentially the property we are interested in,
modulo the following flaw: an $F$-$\invsg$ with a $\zem$
$\theta$ must be a semilattice (Example $1.4.1^{o}$ below);
this is because the maximal element majorizing $\theta$
cannot otherwise be unique. In view of the examples we want
to discuss, we slightly weaken the above requirement, as
follows:

$\ $

{\bf 1.1 Definition} A unital inverse semigroup will be called
an {\em $\widetilde{F}$-$\invsg$} if every $\alpha \in {\cal S}$
which is not zero element for ${\cal S}$ is majorized by a
unique maximal element.

$\ $

The set of maximal elements of an $\widetilde{F}$-$\invsg$ will
be usually denoted by ${\cal M}$.

$\ $

{\bf 1.2 Example} (Toeplitz $\invsg$)  The example which originally
motivated these considerations is the following: let $G$ be a
group, and let $P$ be a unital subsemigroup of $G$. For every
$x \in G$ put:
\begin{equation}
\beta_{x} : \{ t \in P \ | \ xt \in P \} \ni t \
\longrightarrow \ xt \in \{ s \in P \ | \ x^{-1}s \in P \} .
\end{equation}
Clearly, $\beta_{x}$ belongs to the $\siminvsg$ ${\cal I}_{P}$.
Note that $\beta_{x}$ can be the void map $\theta$ on $P$; this
happens if and only if $x \in G \setminus PP^{-1}$.
Let $\sgp$ be the subsemigroup of ${\cal I}_{P}$ generated by
$\{ \beta_{x} \ | \ x \in G \}$; this is a unital
$\star$ -subsemigroup, since
$\beta_{x}^{*} = \beta_{x^{-1}}$, $x \in G$, and since
$\beta_{e}$ ($e$ = unit of $G$) is the identity map on $P$.
Due to its relation to the Toeplitz algebra associated to
$G$ and $P$, we shall call $\sgp$ the Toeplitz $\invsg$ of
$(G,P)$.

Now, if $PP^{-1} = G$, it turns out that $\sgp$ has no $\zem$,
and that it is an $F$-$\invsg$, with
${\cal M} = \{ \beta_{x} \ | \ x \in PP^{-1} \}$.
On the other hand, if $PP^{-1} \neq G$, then $\sgp$ certainly
has a $\zem$ (the void map $\theta$, which is $\beta_{x}$ for
every $x \in G \setminus PP^{-1}$); in this case, $\sgp$
still is an $\f$-$\invsg$, where again
${\cal M} = \{ \beta_{x} \ | \ x \in PP^{-1} \}$
(see Section 3 below).

$\ $

{\bf 1.3 Example} (singly generated $\invsg$) Let ${\cal S}$
be a unital $\invsg$ which is generated (as unital
$\star$-semigroup) by an element $\beta \in {\cal S}$. One
checks immediately, by embedding ${\cal S}$ into a $\siminvsg$,
that we have $\beta^{m} \beta^{n} \leq \beta^{m+n}$ for every
$m,n \in \Z$ (where we make the convention that $\beta ^{n}$
means $\beta^{\star |n|}$ for $n<0$). Since any element of
${\cal S}$ is a product of $\beta$'s and $\beta^{*}$'s, it
immediately follows that for every $\alpha \in {\cal S}$ there
exists (at least one) $n \in \Z$ such that
$\alpha \leq \beta^{n}$.

In many examples (of ${\cal S}$ generated by $\beta$, as above),
the $n \in \Z$ such that $\alpha \leq \beta^{n}$ is uniquely
determined, for every $\alpha \in \Scal$ which is not $\zem$.
This implies that $\Scal$ is an $\widetilde{F}$-$\invsg$, with
$\M = \{ n \in \Z \ | \ \beta^{n}$ is not $\zem$ for $\Scal \}$.

More precisely, let us assume (without loss of generality) that
$\Scal \subseteq {\cal I}_{T}$, the
symmetric inverse semigroup on a non-void set
$T$. Denote the subset $\cap_{n \in \Z}  Dom( \beta^{n} )$ of $T$
by $T_{\infty}$, and put $T_{f} = T \setminus T_{\infty}$. Then
$T_{\infty}$ and $T_{f}$ are invariant for $\beta$, and
$\beta | T_{\infty} $ is a bijection. It is easily seen
that the only situation when
$\Scal$ can fail to be an $\f$-$\invsg$ is when $T_{f}$, $T_{\infty}
\neq \emptyset$, $\beta | T_{\infty}$ is periodic and nonconstant,
and $\beta | T_{f}$ is nilpotent, i.e. there is a power of it
which is the void map. This does not happen, for
instance, in any of the examples discussed in [4, 5].

$\ $

{\bf 1.4 Examples} Various examples of $F$-inverse semigroups
are studied in the algebraic theory of $\invsg$s. For instance:

$1^{o}$ Let ${\cal E}$ be a unital semilattice, i.e. a
unital $\invsg$ all the elements of which are idempotents
(see \cite{Pet}, I.3.9). Every $\alpha \in {\cal E}$ is
selfadjoint (because $\alpha^{3} = \alpha$, hence the unique
$\alpha^{*}$  satisfying
$\alpha \alpha^{*} \alpha = \alpha$,
$\alpha^{*} \alpha  \alpha^{*} = \alpha^{*}$
is $\alpha^{*} = \alpha$). Hence $\alpha \leq \beta$
$\Leftrightarrow$ $\alpha \beta = \alpha$, $\alpha , \beta \in
{\cal E}$, which implies that ${\cal E}$ is an
$\f$-$\invsg$ with $\M = \{ \epsilon \}$.

Let us recall here the basic fact that any two idempotents of
an arbitrary inverse semigroup $\Scal$ are commuting (see \cite{Pet},
II.1.2). This implies that, for every $\invsg$ $\Scal$, the subset
of idempotents $\Szero = \{ \alpha \in \Scal \ | \ \alpha^{2}
= \alpha \}$ is a subsemigroup ( a semilattice).

$2^{o}$ Let $G$ be a group, and let $(G_{n})_{n \in \N}$ be a
sequence of subgroups of $G$, such that $G_{0} =G$ and
$G_{n+1} \subseteq G_{n}$, $n \geq 0$. Consider also
$G_{\infty} = \cap_{n \in \N} G_{n}$. Define $\Scal$ =
$\cup_{0 \leq n \leq \infty} G_{n} \times \{ n \}$, with
multiplication
\[
(x,m)(y,n) \ = \ (xy, \mbox{ min}(m,n)), \ \
0 \leq m,n \leq \infty, \ x \in G_{m}, y \in G_{n} .
\]
This is an example of a Clifford $E$-unitary semigroup (see
\cite{Pet}, II.2 and III.7.1). It is an $F$-$\invsg$, with
$\M =  ( \cup_{n \in \N} (G_{n} \setminus G_{n+1} ) \times
\{ n \} ) \cup ( G_{\infty} \times  \{ \infty \} )$.

$3^{o}$ Let $G$ be a group and let $\sigma :G \rightarrow G$
be an automorphism. Define $\Scal = \N \times G \times \N$,
with multiplication
\[
(m,x,n)(p,y,q) = \left( m-n+ \mbox{max}(n,p),
\sigma^{\mbox{max}(n,p)-n}(x) \sigma^{\mbox{max}(n,p)-p}(y),
q-p+ \mbox{max}(n,p) \right) ,
\]
for $m,n,p,q \in \N$ and $x,y \in G$. This is an example of a
Reilly semigroup (see \cite{Pet}, II.6); it is an
$F$-$\invsg$, with $\M = \{ (m,x,n) \ |$
$m,n \in \N, x \in G, \mbox{ min}(m,n) =0 \}$. In the case when
$G$ has only one element, $\Scal ( \simeq \N \times \N )$ is called
the bicyclic semigroup; this was also appearing in the
context of Example 1.3 (for $(G,P)=( \Z , \N )$, in the notations
used there).

$4^{o}$  It should also be kept in mind that, obviously, every
group is an $F$-$\invsg$ (for an $F$-$\invsg$ $\Scal$ we have
``$\Scal$ group $\Leftrightarrow$ $\M = \Scal$'').

$\ $

{\bf 1.5 The multiplicative structure on $\M$}
Let $\Scal$ be an $\f$-$\invsg$, and let $\M$ be the set of maximal
elements of $\Scal$. It is clear that $\M$ contains the unit
$\epsilon$ of $\Scal$, and that it is closed under $\star$-operation.

We shall work with a multiplicative structure on $\M$, which is
related to the multiplication of $\Scal$, but can't of course coincide
with it. In order to avoid any confusion, we shall use the notation
\begin{equation}
\M \ = \ \{ \beta_{x} \ | \ x \in M \},
\end{equation}
where $M$ is some fixed set of the same cardinality with $\M$
(and $M \ni x \rightarrow \beta_{x} \in \M$ is a bijection), and we
shall define the multiplicative structure we need as an operation
on $M$. We denote by $e$ the unique element of $M$ such that
$\beta_{e} = \epsilon$ (=unit of $\Scal$).

The case when $\Scal$ has not a $\zem$ is quite clear, and well-known
(see \cite{Pet}, III.5, VII.6). We define the multiplication on $M$ by
\begin{equation}
x \cdot y \ = \ \mbox{the unique} \ z \in M \mbox{ such that }
\beta_{x} \beta_{y} \leq \beta_{z}, \ \ x,y \in M.
\end{equation}
Then $M$ becomes a group ($e$ is the unit, and the inverse of
$x \in M$ is the unique $u \in M$ such that
$\beta_{u} = \beta_{x}^{*}$). The map $\Scal \rightarrow M$ which
sends $\alpha \in \Scal$ into the unique $x \in M$ with
$\beta_{x} \geq \alpha$ is, clearly, a surjective semigroup
homomorphism; moreover, it is easily seen that every homomorphism
of $\Scal$ onto a group can be factored by it. Hence $M$ is the
maximal group homomorphic image of $\Scal$.

If $\Scal$ has a $\zem $ $\theta$, then we have on $M$ only a partially
defined multiplication,
\begin{equation}
M^{(2)} = \{ (x,y) \in M \times M \ | \ \beta_{x} \beta_{y}
\neq \theta \} \ \longrightarrow \ M,
\end{equation}
described by the same rule as in (1.4). Still, for all our purposes
this partial multiplication will be as reliable as a group
structure (actually, it is not clear whether it wouldn't be always
possible to embed it into a group). Remark that:

- We still have a partial associativity property; i.e., if we put
\begin{equation}
M^{(3)} = \{ (x,y,z) \in M \times M \times M \ |
\ \beta_{x} \beta_{y} \beta_{z} \neq \theta \} ,
\end{equation}
then $(x \cdot y) \cdot z$ and $x \cdot (y \cdot z)$ make sense and
are equal for every $(x,y,z) \in M^{(3)}$.

- $e$ still is a unit for $M$, i.e. $(x,e), \ (e,x) \in M^{(2)}$
and $e \cdot x = x \cdot e = x$ for all $x \in M$.

- For every $x \in M$, the unique $u \in M$ such that
$\beta_{u} = \beta_{x}^{*}$, which will be denoted by $x^{-1}$,
has $(x^{-1},x), \ (x,x^{-1}) \in M^{(2)}$ and
$x^{-1} \cdot x = x \cdot x^{-1} = e$. Moreover, it is
immediate that $(x^{-1})^{-1} = x$, for every $x \in M$, and
that $(x,y) \in M^{(2)}$ $\Rightarrow$ $(y^{-1},x^{-1}) \in M$,
$y^{-1} \cdot x^{-1} = (x \cdot y)^{-1}$.

- If $(x,y) \in M^{(2)}$ and $x \cdot y =z$, then automatically
$(x^{-1},z), \ (z,y^{-1}) \in M^{(2)}$, and
$x^{-1} \cdot z =y$, $z \cdot y^{-1} =x$. Indeed, we have
$\beta_{x} \beta_{y} \neq \theta$ $\Rightarrow$
$\beta_{x} \beta_{y} \beta_{y}^{*} \neq \theta$ (since
$(\beta_{x} \beta_{y} \beta_{y}^{*}) \beta_{y}$ =
$\beta_{x} \beta_{y}$), hence $(x \cdot y) \cdot y^{-1}$ and
$x \cdot ( y \cdot y^{-1})$ make sense and are equal, by the
partial associativity; but $(x \cdot y) \cdot y^{-1} =
z \cdot y^{-1}$, while $x \cdot (y \cdot y^{-1}) =x$. The
equality $x^{-1} \cdot z =y$ is proved similarly. Note that, as
a consequence, the partial multiplication on $M$ has the
cancellation property ($x \cdot y = x \cdot z$ or
$y \cdot x = z \cdot x$ $\Rightarrow$ $y=z$).

$\ $

$\ $

\section{The groupoid construction} \

{\bf 2.1 Definition} Let $\Scal$ be a unital $\invsg$, and let
$\Omega$ be a locally compact Hausdorff space. By an action
\setcounter{equation}{0}
of $\Scal$ on $\Omega$ we shall understand a $\star$-$\homo$
$\Phi$ of $\Scal$ into the $\siminvsg$ of $\Omega$, such that:

(i) for every $\alpha \in \Scal$, the domain $Dom ( \Phi ( \alpha ))$
and the range $Ran ( \Phi ( \alpha ))$ of $\Phi ( \alpha )$ are
open in $\Omega$, and $\Phi ( \alpha )$ is a $\homeo$ between them;

(ii) $\Phi ( \epsilon ))$ is the identity map on $\Omega$ (where
$\epsilon$ denotes, as usual, the unit of $\Scal$);

(iii) if $\Scal$ has a zero element $\theta$, then $\Phi ( \theta ))$
is the void map on $\Omega$.

$\ $

For the terminology and basic facts about locally compact groupoids
used in this note, we refer the reader to the monograph
\cite{Ren}. The following groupoid construction is a natural
generalization to $\f$-$\invsg$s of the very basic example of groupoid
associated to a group action (\cite{Ren}, Example I.1.2a).

$\ $

{\bf 2.2 Definition} Let $\Scal$ be a unital $\f$-$\invsg$, and let
$\Phi$ be an action of $\Scal$ on the locally compact Hausdorff space
$\Omega$. We use the notations related to $\Scal$ which were
introduced in Section 1.5 above ($M$ such that
$\M = \{ \beta_{x} \ | \ x \in M \}$, and the multiplicative
structure on $M$). We then define a groupoid $\G$, as follows:
\begin{equation}
\G \ = \ \{ (x, \omega ) \ | \ x \in M,
\omega \in Dom ( \Phi ( \beta_{x} )) \subseteq \Omega \} .
\end{equation}
The space of units of $\G$ is $\Omega$, and the domain and range
of $(x, \omega ) \in \G$ are:
\begin{equation}
d(x, \omega ) =  \omega , \
r(x, \omega ) = ( \Phi ( \beta_{x} ))( \omega ) .
\end{equation}
The multiplication on $\G$ is defined by the rule:
\begin{equation}
(x, \omega )(x ', \omega ') \ = \
(x \cdot x ', \omega ' ),
\end{equation}
for $(x, \omega ), \ (x ', \omega ') \in \G$ such that
$( \Phi ( \beta_{x'} ) )( \omega ') = \omega$.
(Note: from the latter equality and
$\omega \in Dom( \Phi ( \beta_{x} ))$
it follows that $\omega '$ is in the domain of
$\Phi ( \beta_{x} ) \Phi ( \beta_{x'} )$ =
$\Phi ( \beta_{x} \beta_{x'} )$; then $\beta_{x} \beta_{x'}$ can't be
zero element for $\Scal$, and using
$\beta_{x} \beta_{x'} \leq \beta_{x \cdot x'}$ we get
$\omega ' \in Dom( \Phi ( \beta_{x \cdot x'} ))$,
$( \Phi ( \beta_{x \cdot x'} ))( \omega ' ) =
( \Phi ( \beta_{x} ))( \omega )$. Thus the right-hand side of
(2.3) is indeed in $\G$, and has
$d(x \cdot x ', \omega ' ) = d(x ', \omega ' )$,
$r(x \cdot x ', \omega ' ) = r(x, \omega )$.)

The identity at the unit $\omega \in \Omega$ is $(e, \omega )$,
with $e$ the unit of $M$, and the inverse of
$(x, \omega ) \in  \G$ is
$(x^{-1}, ( \Phi ( \beta_{x} ))( \omega ) )$.

The topology on $\G$ ($\subseteq M \times \Omega$) is the one
obtained by restricting to $\G$ the product of the discrete
topology on $M$ and the given topology on $\Omega$. This is
locally compact (and Hausdorff), since $\G$ is obviously open
in $M \times \Omega$. It is immediately seen that multiplication
and taking the inverse on $\G$ are continuous, i.e. $\G$ is a
locally compact groupoid. The topology induced from $\G$ to the
space of identities $\Gzero = \{ (e, \omega ) \ |$
$\omega \in \Omega \} \simeq \Omega$ is the one of
$\Omega$. Moreover, $\Gzero$ is open in $\G$, i.e. $\G$ is an
r-discrete groupoid (\cite{Ren}, Definition I.2.6).

$ \ $

{\bf 2.3 Remark}
Note that $\G$ in 2.2 is the disjoint union of the sets
$\{ x \} \times Dom( \Phi ( \beta_{x} ))$, $x \in M$, each
of them open and closed in $\G$. For every $f \in C_{c} ( \G )$,
$supp \ f$ is contained in a finite union
$\cup_{i=1}^{n} \{ x_{i} \} \times Dom( \Phi ( \beta_{x_{i}} ))$,
and we can write $f = \sum_{i=1}^{n} fh_{i}$ with $h_{i}$ the
characteristic function of
$\{ x_{i} \} \times Dom( \Phi ( \beta_{x_{i}} ))$,
$1 \leq i \leq n$ ($fh_{i} \in C_{c}( \G )$, since $h_{i}$ is
continuous). This shows that we have the direct sum decomposition
\begin{equation}
C_{c} ( \G ) \ = \ \oplus_{x \in M} \{ f \in C_{c} ( \G ) \ | \
supp \ f \subseteq \{ x \} \times
Dom( \Phi ( \beta_{x} )) \ \} .
\end{equation}

Note also that from Proposition I.2.8 of \cite{Ren} it follows
that the counting measures form a Haar system on $\G$ (in the
sense of \cite{Ren}, Definition I.2.2).

$\ $

{\bf 2.4 Remark} Let $\Scal$ be a unital $\invsg$ generated by an
element $\beta \in \Scal$, and such that (as in Example 1.3 above),
every $\alpha \in \Scal$ which is not zero element is majorized by
a unique $\beta^{n}, \ n \in \Z$. We have an obvious choice of
notation for the set $M$ appearing in equation (1.3) of 1.5. If
$\beta^{n}$ is not zero element for $\Scal$, for every $n \in \Z$,
then $\Scal$ has no zero element,
$\M = \{ \beta^{n} \ | \ n \in \Z \}$,
and we take $M= \Z$ (the map $M \rightarrow \M$ being
$n \rightarrow \beta^{n}$). Otherwise, $\Scal$ must have a $\zem$
$\theta$, and there exists $m \geq 1$ such that
$\beta^{m} \neq \theta = \beta^{m+1}$,
$\beta^{ {\star} m} \neq \theta = \beta^{ {\star} (m+1)}$;
we then have $\M = \{ \beta^{n} \ | \ |n| \leq m \}$, and we take
$M = \Z \cap [-m,m]$. In any case, the (possibly partially defined)
multiplication on $M$ is the usual addition of integers. For
$\Phi$ an action of $\Scal$ on a space $\Omega$, as in 2.1, note that
the groupoid $\G$ associated to $\Phi$ is
\begin{equation}
\G \ = \ \{ (m, \omega ) \ | \ m \in \Z , \omega \in
Dom( \Phi ( \beta^{n} )) \},
\end{equation}
even in the case when $M$ is of the form $\Z \cap [-m,m]$ (in this
case we have $Dom( \Phi ( \beta^{n} )) = \emptyset$ for
$|n| > m$).

Now, remark that $\G$ of (2.5) is an r-discrete locally compact
groupoid (with the structure defined in 2.2), even if there is no
assumption on $\beta$ to ensure that $\Scal$ is an $\f$-$\invsg$.
Indeed, it is seen immediately that the only things required
for having a valid groupoid structure on $\G$ of (2.5) are the
inequality $\beta^{n} \beta^{m} \leq \beta^{n+m}$,
$n,m \in \Z$, and the group properties of the addition of the
integers, which hold unconditionally.

Moreover, the isomorphism between the C*-algebra of $\G$ in (2.5),
on one hand, and the crossed-product C*-algebra
(in the sense of \cite{E1}) by the partial
automorphism given by $\beta$, on the other hand, will also
turn out to hold with no condition on $\beta$ (see Theorem 4.1
below).

The above considerations suggest an other simple method
of constructing groupoids associated to actions of $\invsg$s.
Let $G$ be a group, let $\Scal$ be a unital inverse
semigroup, and let
$G \ni x \rightarrow \beta_{x} \in \Scal$ be a map such that:
$\beta_{e} = \epsilon$, where $e$ is the unit of $G$ and
$\epsilon$ the one of $\Scal$;
$\beta_{x^{-1}} = \beta_{x}^{*}$ for every $x \in G$;
$\beta_{x} \beta_{y} \leq \beta_{xy} $ for every $x,y \in G$;
and for every $\alpha \in \Scal$ there exists (at least one)
$x \in G$ such that $\alpha \leq \beta_{x}$.
(In other words, instead of giving conditions which ensure a
multiplicative structure on a remarkable family of elements of
$\Scal$, we impose this from outside.) Then, to an action $\Phi$
of $\Scal$ on a locally compact space $\Omega$, one can associate
the locally compact groupoid
\begin{equation}
\G \ = \ \{ (x, \omega ) \ | \ x \in G, \omega \in
Dom( \Phi ( \beta_{x} )) \} ,
\end{equation}
with groupoid structure defined as in 2.2. Note that the actions
of Toeplitz inverse semigroups can also be considered in this
way (in the notations of Example 1.2, we take
$G \ni  x \rightarrow \beta_{x} \in \sgp$ to be exactly the one
given by equation (1.2)); however, the properties shown in
Section 6 below don't hold in general for groupoids of the type
(2.6).

$\ $

\section{Example: the Toeplitz inverse semigroup}

We now return to the framework of Example 1.3. Consider
\setcounter{equation}{0}
$G, P, (\beta_{x} )_{x \in G}$ and
$\sgp \subseteq {\cal I}_{P}$ as in 1.3.
We begin by proving the assertions which were made there about $\sgp$.

$\ $

{\bf 3.1 Lemma} If $\pq = G$, then $\sgp$ has no zero element.

$\ $

{\bf Proof} If $\sgp$ has a zero element $\theta$, then this must be the void
map on $P$ (we leave here apart the case when $G = P = \{ e \}$).
Indeed, $\theta$ would otherwise be the identical transformation of a subset of
$P$ which is fixed by every $\beta_{x}$, $x \in \pq$;
but for $x \neq e$, $\beta_{x}$ has no fixed points.

Hence, since an arbitrary element of $\sgp$ is a product of $\beta_{x}$'s,
what we need to show is that for every $n \geq 1$ and $\xuri \in G$, the
partially defined transformation $\beturi$ on $P$ is non-void.

Let $\xuri \in G$ be arbitrary. Since $\pq = G$, we can find
$s_{1}, \ldots , s_{n}, t_{1}, \ldots , t_{n} \in P$ such that:
$x_{n} = s_{n} t_{n}^{-1}$, $x_{n-1} s_{n} = s_{n-1}t_{n-1}^{-1}$,
$x_{n-2} s_{n-1} = s_{n-2}t_{n-2}^{-1}$, $\ldots$,
$x_{1}s_{2} = s_{1}t_{1}^{-1}$.
Then $t = t_{n} \dots t_{1} \in P$ is in the domain of
$\beturi$, because, as it is easily checked,
$x_{j} \dots x_{n} t = s_{j} t_{j-1} \dots t_{1} \in P$, $1 \leq j \leq n$.
Thus $\beturi$ is indeed non-void. {\bf QED}

$\ $

{\bf 3.2 Lemma} For every non-void $\alpha \in \sgp$ there exists
a unique $x \in \pq $ such that $\alpha \leq \beta_{x}$; this $x$ can
be expressed as $\alpha (t) t^{-1}$, with $t$ arbitrary in the domain of
$\alpha$.

$\ $

{\bf Proof} Write $\alpha = \beturi$, with $\xprod \in \pq$; then
$Dom( \alpha ) = \{ t \in P \ |$
$x_{n} t , x_{n-1} x_{n} t , \ldots , \xprod t \in P\}$,
and $\alpha (t) = \xprod t $ for $t \in Dom(\alpha )$.
Putting $x = \xprod $, we get $\alpha (t) t^{-1} = x$ for every
$t \in Dom (\alpha )$ (in particular $x$ is in $\pq$). It is clear that
$\alpha \leq \beta_{x}$, and that $x$ is the unique element of
$\pq$ having this property. {\bf QED}

$\ $

{}From the above two lemmas it is clear that no matter if $\pq =G$ or
not, $\sgp$ is an $\f$-$\invsg$, with
$\M = \{ \beta_{x} \ | \ x \in \pq \}$.
It fits very well the notations of Section 1.5 to take
$M = \pq$. Note also that the multiplication defined on $M = \pq$ by
equation (1.4) of 1.5 coincides with the one coming from the group $G$;
this happens because for any $x, y \in \pq$ such that
$\beta_{x} \beta_{y} \neq \theta$, $\beta_{x} \beta_{y}$ acts on
its domain by left translation with $xy$.

$\ $

We now pass to describe a remarkable action of $\sgp$, which gives
via the construction of 2.2 the Wiener-Hopf
groupoid $\G$, introduced by
P.Muhly and J.Renault in \cite{MR}.
The space of the action of $\sgp$ will be
\begin{equation}
\Omega  \ = \ \ clos  \{ t P^{-1} \ | \
t \in P \} \subseteq \{ 0,1 \}^{G} ,
\end{equation}
where $\{ 0,1 \}^{G}$ is identified
to the space of all subsets of $G$. The
topology on $\Omega$ is the restriction of the product topology
on $\{ 0,1 \}^{G}$, and is compact and Hausdorff. Note that
$P^{-1} \subseteq A \subseteq \pq$ for every $A \in \Omega$.

$\ $

{\bf 3.3 Proposition} For every $\alpha \in \sgp $ we define:
\begin{equation}
\left\{ \begin{array}{l}
\Phi (\alpha ) :
\ clos  \{t P^{-1} \ | \ t \in Dom(\alpha ) \} \rightarrow
\ clos \{s P^{-1} \ | \ s \in Ran(\alpha ) \}, \\
(\Phi (\alpha )) (A) = xA,
\end{array}
\right.
\end{equation}
where $x$ in (3.2) is the unique element of $\pq$ such that
$\alpha \leq \beta_{x}$. (If $\alpha = \theta $, the void map
on $P$, we take
by convention $\Phi (\alpha )$ to be the void map on $\Omega$.) Then
$\Phi (\alpha )$ makes sense for every $\alpha \in \sgp$,
and $\Phi$ is an action of
$\sgp$ on $\Omega$, in the sense of Definition 2.1.

$\ $

{\bf Proof} Let $\alpha \in \sgp$ be non-void,
and let $x$ be the unique
element of $\pq$ such that $\alpha \leq \beta_{x}$. The map
$A \rightarrow xA$ (=$\{ xa \ | \ x \in A \}$) is continuous on
${\{ 0,1 \} }^{G}$, hence
$\{ A \subseteq G \ | \ xA \in$ $clos \{ sP^{-1} \ |$
$s \in Ran( \alpha ) \} \ \}$ is closed. This set contains $tP^{-1}$
for every $t \in Dom( \alpha )$ (because $t \in Dom ( \alpha )
\Rightarrow xt = \alpha (t) \in Ran( \alpha )$, and on the other hand
$x(tP^{-1})=(xt)P^{-1}$); so $\Phi ( \alpha )$ defined in (3.2)
makes sense. It is also clear that $\Phi ( \alpha^{*} )$ is an
inverse for $\Phi ( \alpha )$, which is thus a bijection.

It is useful to note that if $\alpha \in \sgp$ is written as a
product $\beturi$ (for some $n \geq 1$ and $\xuri \in \pq$), then
we have the equivalent characterization
\begin{equation}
Dom( \Phi ( \alpha )) \ = \ \{ A \in \Omega \ | \
A^{-1} \ni x_{n}, \ x_{n-1}x_{n}, \ \ldots \ , \xprod \}
\end{equation}
(note: this is valid even if
$\alpha = \theta $). Indeed, the right-hand
side of (3.3) is immediately seen to be equal to $clos \{ tP^{-1} \ |$
$(tP^{-1})^{-1} \ni x_{n},x_{n-1}x_{n} , \ldots , \xprod \}$ =
$clos \{ tP^{-1} \ |$
$t, x_{n}t, x_{n-1}x_{n}t, \ldots , \xprod t \in P \}$,
which is exactly $Dom( \Phi ( \alpha ))$.

As a consequence of (3.3), it is clear that $Dom ( \Phi ( \alpha ))$
(and $Ran( \Phi ( \alpha )) = Dom ( \Phi ( \alpha^{*} ))$, too) is
open in $\Omega$ for every $\alpha \in \sgp$.

We are left to show that
$\Phi ( \alpha ) \Phi ( \alpha ' ) =
\Phi ( \alpha \alpha ' )$ for every $\alpha , \alpha ' \in \sgp$.
If one of $\alpha , \alpha '$ is $ \theta$, then both
$\Phi ( \alpha ) \Phi ( \alpha ' )$ and
$\Phi ( \alpha \alpha ' )$
are the void map on $\Omega$, so we shall assume
$\alpha \neq \theta \neq \alpha '$ (but we don't assume
$\alpha \alpha ' \neq \theta$). We take
$x_{1}, \ldots ,x_{m}, y_{1}, \ldots , y_{n} \in \pq$ such that
$\alpha = \beta_{x_{1}} \cdots \beta_{x_{m}}$,
$\alpha ' = \beta_{y_{1}} \cdots \beta_{y_{n}}$; note that the
unique $x,y \in \pq$ such that $\alpha \leq \beta_{x}$,
$\alpha ' \leq \beta_{y}$ must then be
$x = x_{1} \cdots x_{m}$ and $y = y_{1} \cdots y_{n}$. Using
(3.3) we see that
\[
Dom( \Phi ( \alpha ) \Phi ( \alpha ' ) ) \ = \
\left\{ A \in \Omega \ | \
\begin{array}[t]{l}
A^{-1} \ni y_{n}, \ y_{n-1}y_{n}, \ \ldots \ , y_{1} \cdots y_{n}, \\
(yA)^{-1} \ni x_{m}, \ x_{m-1}x_{m}, \ \ldots \ , x_{1} \cdots x_{m}
\end{array}
\right\}
\]
\[
= \ \left\{ A \in \Omega \ | \
\begin{array}[t]{l}
A^{-1} \ni y_{n}, \ y_{n-1}y_{n}, \ \ldots \ , y_{1} \cdots y_{n}, \\
x_{m} \igreci , \ x_{m-1}x_{m} \igreci , \
\ldots \ , x_{1} \cdots x_{m} \igreci
\end{array}
\right\}
\]
\[
= \ Dom( \Phi ( \alpha \alpha ' ))
\]
(where the latter equality holds also by (3.3), since
$\alpha \alpha '$ =
$\beta_{x_{1}} \cdots \beta_{x_{m}}
\beta_{y_{1}}  \cdots \beta_{y_{n}}$).
If $Dom ( \Phi ($ $\alpha )  \Phi ( \alpha ' ))$ =
$Dom ( \Phi ( \alpha \alpha ' ))$ is non-void, it is clear that
both transformations act on it
by left translation with $xy$, hence they are equal.
{\bf QED}

$\ $

Note that if $P \cap P^{-1} = \{ e \}$, then $\Omega$ of (3.1) is a
compactification of $P$ (by identifying $t \in P$ with
$tP^{-1}$); in some sense, $\Omega$ is obtained from $P$ by adding
one point for each ``type of convergence to $\infty$ in $P$''
(compare to the comments in Section 2.3.1 of \cite{N1}). Taking
into account that a non-void $\alpha$ in $\sgp$ is actually
the left translation with $x$ on $Dom( \alpha )$, with $x \in \pq$
such that $\alpha \leq \beta_{x}$, we can interprete the map
$\Phi ( \alpha ) $ of (3.2) as a sort of ``compactification of
$\alpha$''.
\footnote[1]{ One may ask what happens if we don't do any
compactification, and just take the obvious action of $\sgp$
on $P$. It is immediate that the groupoid associated to this
would be the total equivalence relation on $P$, having thus
C*-algebra isomorphic to the compact operators on $l^{2}(P)$. }

$\  $

Finally (without having to assume $P \cap P^{-1} = \{ e \}$), we
have, in view of the characterization (3.3):
\begin{equation}
Dom( \Phi ( \beta_{x} )) \ = \
\{ A \in \Omega \ | \ A \ni x^{-1} \} , \ \  x \in \pq .
\end{equation}
Hence the groupoid associated as in 2.2 to the above action
of $\sgp$ on $\Omega$ is:
\begin{equation}
\G \ = \ \{ (x,A) \ | \ A \in \Omega , x \in A^{-1}  \} .
\end{equation}
This is exactly the form given to the Wiener-Hopf groupoid
in \cite{N1}, Proposition 2.3.4.

$\ $

\section{Example: singly generated inverse semigroups} \

Let $\Scal$ be a unital $\invsg$, generated (as unital
\setcounter{equation}{0}
$\star$-semigroup) by an element $\beta \in \Scal$, and let
$\Phi$ be an action of $\Scal$ on the locally compact space
$\Omega$ (as in 2.1). We denote $\Phi ( \beta )$ by
$\tilbeta$. Let $\G = \{ (n, \omega ) \ |$
$n \in \Z , \ \omega \in Dom( \tilbeta^{n} ) \}$ be the groupoid
associated to this action, as in equality (2.5) of
Remark 2.4.

On the other hand, if we consider the ideals
$I = \{ f \in C_{o} ( \Omega ) \ | \ f \equiv 0$ on
$\Omega \setminus Ran( \tilbeta ) \}$,
$J = \{ f \in C_{o} ( \Omega ) \ | \ f \equiv 0$ on
$\Omega \setminus Dom( \tilbeta ) \}$ of $C_{o} ( \Omega )$,
then $\tilbeta$ determines an isomorphism
$\theta : I \rightarrow J$,
\[
( \theta (f))( \omega ) \ = \ \left\{
\begin{array}{ll}
f( \tilbeta ( \omega )),
& \mbox{  if } \omega \in Dom( \tilbeta )  \\
0,                       & \mbox{  otherwise}
\end{array}  \right.
\ \ , \ \mbox{for } f \in I .
\]
Thus $\Theta = ( \theta , I,J)$ is a partial automorphism
of $C_{o} ( \Omega )$, in the sense of R. Exel \cite{E1},
Definition 3.1, and has attached to it a covariance
C*-algebra C*($C_{o} ( \Omega ), \Theta )$ (see \cite{E1},
Definition 3.7).

$\ $

{\bf 4.1 Theorem} Assuming $\Omega$ second countable, we have
that C*($C_{o}, \Theta )$ and C*($\G$) are naturally isomorphic.

$\ $

{\bf Proof} Following the notations of \cite{E1}, Section 3,
let us put
\[
D_{n} \ = \ \{
f \in C_{o} ( \Omega ) \ | \ f \equiv 0 \mbox{ on }
\Omega \setminus Dom( \tilbeta^{n} ) \} , \ \  n \in \Z ;
\]
in other words, $D_{n}$ is the domain of $\theta^{-n}$. We have
in particular $D_{o} = C_{o} ( \Omega )$, $D_{1} =J$, $D_{-1} =I$.
Let $L$ be the vector space of sequences $(f_{n})_{n \in \Z}$,
with $f_{n} \in D_{n}$ for every $n$, and such that
$\sum_{n \in \Z} ||f_{n} || < \infty$; for $m \in \Z$ and
$f \in D_{m}$ denote by $f \delta_{m}$ the sequence
$(f_{n})_{n \in \Z}$
in $L$ which has $f_{m} = f$ and $f_{n} = 0$ for $ n \neq m$.
Then $L$ is given (in \cite{E1}, p.7) a $\star$-algebra
structure which in particular has for $n,m \in \Z$,
$f_{n} \in D_{n}$, $f_{m} \in D_{m}$:
$(f_{n} \delta_{n}) \star (f_{m} \delta_{m})$ =
$g \delta_{n+m}$,
$(f_{n} \delta_{n})^{*} = h \delta_{-n}$, with:
\begin{equation}
g( \omega ) \ = \ \left\{  \begin{array}{ll}
f_{n} ( \omega ) f_{m} ( \tilbeta^{n} ( \omega )),  &
\mbox{  if } \omega \in Dom( \tilbeta^{n} ) \cap
Dom( \tilbeta^{n+m} )   \\
0,   & \mbox{  otherwise,}
\end{array}  \right.
\end{equation}
\begin{equation}
h( \omega ) \ = \ \left\{  \begin{array}{ll}
\overline{f_{n} ( \tilbeta^{-n} ( \omega ))},  &
\mbox{  if } \omega \in Dom( \tilbeta^{-n} )    \\
0,   & \mbox{  otherwise.}
\end{array}  \right.
\end{equation}
C*($C_{o} ( \Omega ), \Theta )$
is defined as the enveloping C*-algebra of $L$,
with respect to the norm
$|| (f_{n})_{n \in \Z} ||_{1}$ =
$\sum_{n \in \Z} ||f_{n} ||$.

For every $n \in \Z$ let us define
\[
\dnciup \ = \ \{
f \in C_{c} ( \Omega ) \ | \ supp \ f
\subseteq Dom( \tilbeta^{n} ) \}
\subseteq D_{n} ,
\]
and let $\lciup$ be the space of finitely supported
sequences $(f_{n})_{n \in \Z}$ with the property that
$f_{n} \in \dnciup$ for every $n$. Then clearly
$\lciup$ is a $\star$-subalgebra of $L$, dense in
$|| \cdot ||_{1}$, and
C*($C_{o} ( \Omega ), \Theta )$
can also be defined as the enveloping C*-algebra of
$( \lciup , || \cdot ||_{1} )$. Since $\Omega$ is
assumed to be second countable, the latter normed
$\star$-algebra is separable, and we may consider only
its representations on separable Hilbert spaces.

Now, every $\star$-representation
$\pi : \lciup \rightarrow {\cal L} ( {\cal H} )$
is automatically contractive with respect to
$|| \cdot ||_{1}$. In order to verify this, it clearly
suffices to check that
$|| \pi (f_{n} \delta_{n} ) || \leq ||f_{n} ||$
for every $n \in \Z$ and $f_{n} \in \dnciup$. And
indeed, one sees immediately (from (4.1),(4.2)) that
$(f_{n} \delta_{n}) \star (f_{n} \delta_{n})^{*}$ =
$|f_{n}|^{2} \delta_{o}$, hence
$|| \pi (f_{n} \delta_{n} )||^{2}$ =
$|| \pi ( |f_{n} |^{2} \delta_{o} ) ||$; the latter
quantity does not exceed
$|| \ |f_{n}|^{2} \ || = ||f_{n}||^{2}$, because
$\pi$ restricted to
$\{ f \delta_{o} \ | \ f \in C_{c} ( \Omega ) \}
\subseteq \lciup$ gives a $\star$-representation
of $C_{c} ( \Omega )$, which is automatically
contractive. (The latter assertion holds even if
$\Omega$ is non-compact, due to the fact that, for
every $f \in C_{c} ( \Omega )$, the spectral radius
of $f$ in the unitization of $C_{c} ( \Omega )$
equals $||f||$. )

It results that C*($C_{o} ( \Omega ), \Theta )$
is the enveloping C*-algebra of $\lciup$, considered
with respect to all the algebraic $\star$-representations
of $\lciup$ on separable Hilbert spaces.

On the other hand, due to the separability condition
imposed on $\Omega$, it is clear that $\G$ is second
countable; hence, by Corollary II.1.22 of \cite{Ren}
(see also Corollaire 4.8 of \cite{R}), C*($\G$) is
the enveloping C*-algebra of $C_{c}( \G )$ with respect
to all its (algebraic) $\star$-representations on
separable Hilbert spaces. This makes clear that the
Theorem will follow as soon as we can prove that the
$\star$-algebras $\lciup$ and $C_{c}( \G )$ are
isomorphic.

For every $n \in \Z$ and $f_{n} \in \dnciup$ let us
denote by $\chi_{n} \otimes f_{n}$ the restriction to
$\G \subseteq \Z \times \Omega$ of the direct product
between $\chi_{n}$ = (characteristic function of
$\{ n \}$) and $f_{n}$. We denote, for every $n \in \Z$,
${\cal D}_{n}$ =
$\{ \chi_{n} \otimes f_{n} \ |$
$f_{n} \in \dnciup \} \subseteq C_{c} ( \G )$. From the
direct sum decomposition of equation (2.4) it is
immediate that $C_{c}( \G ) \ \simeq$
$\oplus_{n \in \Z} {\cal D}_{n}$. (Note that some of the
spaces ${\cal D}_{n}$ may be reduced to $\{ 0 \}$, if
$\tilbeta$ is nilpotent; in this case, the corresponding
$D_{n}$'s are also $\{ 0 \}$.) In view of the obvious
decomposition $\lciup$ =
$\oplus_{n \in \Z} \{ f_{n} \delta_{n} \ |$
$f_{n} \in \dnciup \}$, it becomes clear that $\lciup$
and $C_{c} ( \G )$ are naturally isomorphic as vector
spaces.

Recalling the definition of the multiplication and
$\star$-operation on $C_{c} ( \G )$ (from \cite{Ren},
Proposition II.1.1) one gets, for $m,n \in \Z$,
$f_{n} \in \dnciup$, $f_{m} \in \widehat{D_{m}}$, the
formulae: $( \chi_{n} \otimes f_{n} )  \star
( \chi_{m} \otimes f_{m} )$ = $\chi_{n+m} \otimes g$,
$( \chi_{n} \otimes f_{n} )^{*} = \chi_{-n} \otimes h$,
where:
\begin{equation}
g( \omega ) \ = \ \left\{  \begin{array}{ll}
f_{n} ( \tilbeta^{m} ( \omega )) f_{m} ( \omega ), &
\mbox{  if } \omega \in Dom( \tilbeta^{m} ) \cap
Dom( \tilbeta^{n+m} )   \\
0,   & \mbox{  otherwise,}
\end{array}  \right.
\end{equation}
and $h$ is the same as in (4.2).

Taking into account
the difference of choice which appears in the
definition of the multiplication in the two approaches
(i.e. in (4.1) vs (4.3)), one has thus to proceed
as follows: for every $n \in \Z$ and $f_{n} \in \dnciup$
denote by $\Gamma_{n} (f_{n})$ the complex conjugate of
the function in (4.2). Then the linear isomorphism
$\lciup \rightarrow C_{c} ( \G )$ determined by
$f_{n} \delta_{n} \rightarrow \chi_{-n} \otimes
\Gamma_{n} (f_{n})$, $n \in \Z$, $f_{n} \in \dnciup$,
is also an isomorphism of $\star$-algebras (the easy
verification of this is left to the reader).
{\bf QED}

$\ $

{\bf 4.2 Remark} Crossed-products by
partial isomorphisms were used in \cite{E2} to approach
AF-algebras, and in particular to approach in a very
explicit way UHF-algebras. This leads to an interesting
point of view on the Glimm groupoid (and more generally
on AF-groupoids, defined in \cite{Ren}, Section III.1),
as being close, in some sense, to transformation groups.

More precisely, let $(n_{i})_{i=0}^{\infty}$ be a sequence
of positive integers, and let
$X = \prod_{i=0}^{\infty} \{ 0,1, \ldots$,
 $n_{i}$ $-1 \}$ have
the product topology. Consider, as in \cite{E2}, the
``restricted odometer map'' $\beta^{*}$ :
$X \setminus \{ (n_{1}-1,n_{2}-1, \ldots  ) \}$
$\rightarrow$ $X \setminus \{ (0,0, \ldots  ) \}$
which is given by the formal addition with
$(1,0,0, \ldots )$, with carry over to the right.
Let $\beta$ be the inverse of $\beta^{*}$, let $\Scal$ be the
unital $\star$-semigroup generated by $\beta$ in
${\cal I}_{X}$, and let $\G$ be the groupoid associated
to the action $\Phi ( \alpha ) \equiv \alpha$ of $\Scal$ on $X$.
Since $\beta$ is very close to be a $\homeo$ of $X$,
$\G$ is in some sense close to be a transformation group.
On the other hand, $\G$ is easily seen to be isomorphic
to the Glimm groupoid on $(n_{i})_{i=0}^{\infty}$ (defined
in \cite{Ren}, p.128); the fact that C*($\G$) $\simeq$
UHF($(n_{i})_{i=0}^{\infty}$) can be verified either this
way, or by combining the above Theorem 4.1 with Theorem 3.2
of \cite{E2}.

$\ $

\section{The relation with localizations (in the sense of
Kumjian)} \

A. Kumjian (\cite{K}) has developed a method which associates a
\setcounter{equation}{0}
C*-algebra to an inverse semigroup $\Scal$ of
$\homeo$s between open subsets of a locally compact space
$\Omega$, such that:
\begin{equation}
(Dom( \alpha ) )_{\alpha \in \Scal}
\mbox{ is a basis for the topology of } \Omega .
\end{equation}
Following \cite{K} Definition 2.3,
such an $\invsg$ will be called a localization.

The point of view of this note is slightly different from the
one of \cite{K}; for instance, the approach to the Glimm groupoid
mentioned in Remark 4.2 above is different from the one taken in
\cite{K} (see 5.2 below); actually, the $\invsg$ generated by the
restricted odometer map is not a $\local$, and so is also
the case for the example described in Section 3 (even in the
classical situation when $(G,P) = (\Z , \N )$ ).

Still, the C*-algebra construction of \cite{K} coincides with
the C*-algebra of the groupoid defined in 2.2, on the common
territory of the two approaches:

$\ $

{\bf 5.1 Theorem} Let $\Scal$ be a countable $\f$-$\invsg$ of
$\homeo$s between open subsets of the locally compact space
$\Omega$, and assume that $\Scal$ is a $\local$.
Let $\G$ be the groupoid associated (as in 2.2) to the
action $\Phi ( \alpha ) \equiv \alpha$ of $\Scal$ on $\Omega$.
Then C*($\G$) is isomorphic to the C*-algebra associated to
$\Scal $ in \cite{K}.

$\ $

{\bf Proof} The C*-algebra associated to $\Scal$ in
\footnote[1]{ We warn the reader about the unfortunate
coincidence that $\Omega$ is used in the present paper to
denote a space, and in \cite{K} to denote an $\invsg$. }
\cite{K} is the envelopation of a $\star$-algebra which we
will denote (following \cite{K}, Section 5) by
$C_{c} ( \Scal )$. An argument similar to the one used in the
proof of Theorem 4.1 shows that it suffices to prove the
isomorphism of $\star$-algebras
$C_{c} ( \Scal) \simeq C_{c} ( \G )$.

We now have to go into the details of the definition of
$C_{ c} ( \Scal )$. Following \cite{K}, Notation 5.2, let us put
\begin{equation}
D( \Scal ) \ = \ \oplus_{\alpha \in \Scal }  \{
( \alpha , f) \ | \ f \in C_{c} (Dom ( \alpha )) \},
\end{equation}
where $C_{c}(Dom ( \alpha ))$ stands for
$\{ f \in C_{c} ( \Omega ) \ |$
$supp \ f \subseteq Dom( \alpha ) \}$, and where the summand in (5.2)
corresponding to $\alpha$ is given the linear structure coming
from $C_{c}(Dom ( \alpha ))$. For $\alpha \in \Scal$ and
$f \in C_{c}(Dom ( \alpha ))$ we shall view $( \alpha ,f)$ as
an element of $D( \Scal )$, in the obvious manner. $D( \Scal)$
is given (in \cite{K}, 5.2) a $\star$-algebra structure, such
that for $\alpha_{1}, \alpha_{2} \in \Scal$,
$f_{1} \in C_{c}(Dom ( \alpha_{1} ))$,
$f_{2} \in C_{c}(Dom ( \alpha_{2} ))$ we have
$(\alpha_{1},f_{1})(\alpha_{2},f_{2})=(\alpha_{1}\alpha_{2},g)$,
$(\alpha_{1},f_{1})^{*} = (\alpha_{1}^{*} ,h)$, with:
\begin{equation}
g( \omega ) \ = \ \left\{  \begin{array}{ll}
f_{1}( \alpha_{2} ( \omega ))f_{2} ( \omega ),  &
\mbox{ if } \omega \in Dom( \alpha_{2} )  \\
0,    & \mbox{ otherwise },
\end{array}  \right.
\end{equation}
\begin{equation}
h( \omega ) \ = \ \left\{  \begin{array}{ll}
\overline{f_{1}( \alpha_{1}^{*} ( \omega ))},  &
\mbox{ if } \omega \in Dom( \alpha_{1}^{*} )  \\
0,    & \mbox{ otherwise }.
\end{array}  \right.
\end{equation}
Then $C_{c} ( \Scal )$ is $D( \Scal )/I_{o}( \Scal )$
(\cite{K}, p.160), where the definition of the ideal
$I_{o} ( \Scal ) \subseteq D( \Scal )$ remains to be
recalled.

On the other hand, let us also consider the groupoid
$\G = \{ (x, \omega ) \ |$
$\ x \in M, \omega \in Dom( \beta_{x}) \}$ defined in 2.2,
where $\{ \beta_{x} \ | \ x \in M \}$ is the set of
maximal elements of $\Scal$, as in 1.5. For every
$x \in M$ and $f \in C_{c} ( Dom(\beta_{x} ))$ let us
denote by $\chi_{x} \otimes f \in C_{c} ( \G )$ the
restriction to $\G ( \subseteq M \times \Omega )$ of
the direct product between $\chi_{x}$ = (characteristic
function of $x$) and $f$.

For $\alpha \in \Scal$ which is not $\zem$, consider the
unique $x \in M$ such that $\alpha \leq \beta_{x}$; then
$Dom ( \alpha ) \subseteq Dom( \beta_{x} )$, hence we have
a linear map $(\alpha ,f) \rightarrow \chi_{x} \otimes f$ from
$\{ ( \alpha ,f) \ | \ f \in C_{c}(Dom( \alpha )) \}$ into
$C_{c}( \G )$. The direct sum (after $\alpha \in \Scal$)
of all these linear maps is a linear map
$D( \Scal ) \rightarrow C_{c}( \G )$, which will be denoted
by $\rho$. Comparing (5.3),(5.4) with the definition of the
operations on $C_{c}( \G )$ (given by formulae similar to
(4.3),(4.2) of Section 4), one checks immediately that
$\rho$ is actually a $\star$-algebra $\homo$. $\rho$ is
clearly onto, and we are left to show that $Ker \ \rho$
equals $I_{o}( \Scal )$, the ideal of $D( \Scal )$
mentioned following to (5.4).

The definition of $I_{o}( \Scal )$ depends on the following
notion: a finite family $( \alpha_{a}, f_{a})_{a \in A}$
(with $\alpha_{a} \in \Scal$,
$f_{a} \in C_{c}(Dom( \alpha_{a} ))$
for $a \in A$) is called coherent if there are open subsets
$(U_{a})_{a \in A}$ of $\Omega$ such that:
\begin{equation}
\left\{  \begin{array}{l}
supp \ f_{a} \subseteq U_{a} \subseteq Dom(\alpha_{a}),
\ \ a \in A, \\
\alpha_{a} | U_{a} \cap U_{a'} =
\alpha_{a'} | U_{a} \cap U_{a'} \mbox{ for } a, a' \in A
\mbox{ having } U_{a} \cap U_{a'} \neq \emptyset .
\end{array}  \right.
\end{equation}
Using this notion, one arrives to $I_{o} ( \Scal )$ in
several steps:

- Define $I( \Scal )$ to be the linear span of
$\{ \sum_{a \in A} ( \alpha_{a},f_{a} ) \ |$
$(\alpha_{a},f_{a})$ coherent, $\sum_{a \in A} f_{a} =0 \}$;
$I( \Scal )$ is shown to be a two-sided, selfadjoint ideal
of $D( \Scal )$ (\cite{K}, p. 158).

- For every $\xi \in D( \Scal )$ consider the number
\[
| \xi |_{o}^{'} \ = \ \mbox{inf } \sum_{b \in B}  ||
\sum_{a \in A_{b}} f_{a,b} ||_{\infty} ,
\]
where the infimum is taken after all the decompositions
$\xi =  \sum_{b \in B} \sum_{a \in A_{b}} ( \alpha_{a,b},
f_{a,b} )$ (with $B$ and $(A_{b})_{b \in B}$ finite sets),
such that each of the families
$( \alpha_{a,b},f_{a,b})_{a \in A_{b}}$
($b \in B$) is coherent.

- For every $\xi \in D( \Scal )$ consider the number
\[
| \xi |_{o} \ = \ \mbox{inf }_{\eta \in I( \Scal )}
| \xi - \eta |_{o}^{'};
\]
then $| \cdot |$ is a $\star$-algebra seminorm on
$D( \Scal )$ (\cite{K}, p.159).

- Define $I_{o} ( \Scal ) = \{ \xi \in D( \Scal ) \ ;
\ | \xi |_{o} = 0 \}$.

$\ $

An obvious situation
of coherent family $( \alpha_{a},f_{a})_{a \in A}$ can be
obtained by taking all the $\alpha_{a}$'s ($a \in A$) to be
majorized by the same $\beta_{x}, \ x \in M$; we shall call
such a family majorized-coherent.
It is immediate that an element $\xi \in D( \Scal )$ belongs
to $Ker \  \rho$ if and only if it is of the form
$\sum_{i=1}^{n} \sum_{a \in A_{i}} ( \alpha_{a}, f_{a} )$,
where each of the families
$(\alpha_{a},f_{a})_{a \in A_{i}}$ ($ 1 \leq i \leq n$)
is majorized-coherent. This makes clear that
$Ker \ \rho \subseteq I( \Scal )$; the opposite inclusion
is also true, as implied by the following

$\ $

{\bf 5.1.1 Lemma} If $( \alpha_{a},f_{a})_{a \in A}$ is a
coherent family, with $\sum_{a \in A} f_{a} =0$, then there
exists a partition $A = A_{1} \cup \ldots \cup A_{n}$ of $A$
such that each of the families
$(\alpha_{a},f_{a})_{a \in A_{j}}$ ($1 \leq j \leq n$) is
majorized-coherent with $\sum_{a \in A_{j}} f_{a} =0$.

$\ $

{\bf Proof of Lemma 5.1.1} Without loss of generality, we
may assume that $\alpha_{a}$ is non-zero for every $a \in A$;
then, for every $a \in A$, we can consider the unique
$x(a) \in M$ such that $\alpha_{a} \leq \beta_{x(a)}$, and
we can partition $A$ such that $a, a'$ are in the same
component of the partition if and only if $x(a)=x(a')$.

Pick a family $(U_{a})_{a \in A}$ of open subsets
of $\Omega$ such that (5.5) holds. We will show that
$x(a) \neq x(a') \Rightarrow U_{a} \cap U_{a'} = \emptyset$;
this immediately implies that the partition considered in
the preceding paragraph satisfies the required conditions.

So, let $a,a'$ be in $A$ such that
$U_{a} \cap U_{a'} \neq \emptyset$.
Using (5.1), we can find a non-zero idempotent $\gamma \in \Scal$
such that $Dom( \gamma ) \subseteq U_{a} \cap U_{a'}$.
Then using (5.5) we get
$\alpha_{a} \gamma = \alpha_{a'} \gamma$, non-zero. We have
$\alpha_{a} \gamma \leq \beta_{x(a)}$,
$\alpha_{a'} \gamma \leq \beta_{x(a')}$; but the maximal
element of $\Scal$ majorizing
$\alpha_{a} \gamma = \alpha_{a'} \gamma$  is unique, hence
$x(a)=x(a')$.

$\ $

We have thus obtained $I( \Scal ) = Ker \ \rho$. It is
obvious that $\xi \in I( \Scal )$ $ \Rightarrow$
$| \xi |_{o}^{'} =0$ $\Rightarrow$ $| \xi |_{o} =0$,
hence $I_{o} ( \Scal ) \supseteq I( \Scal )$. On the other
hand, the inclusion $I_{o} ( \Scal ) \subseteq Ker \ \rho$
will come out from the following

$\ $

{\bf 5.1.2 Lemma:} If $\xi = \sum_{b \in B}
\sum_{a \in A_{b}} ( \alpha_{a,b},f_{a,b}) \in D ( \Scal )$
with $( \alpha_{a,b},f_{a,b})_{a \in A_{b}}$ coherent for
every $b \in B$, and if $x \in M$,
$\omega \in Dom( \beta_{x} )$, then
\begin{equation}
|( \rho ( \xi ))(x, \omega )| \ \leq \
\sum_{b \in B} | \sum_{a \in A_{b}} f_{a,b} ( \omega ) | .
\end{equation}

$\ $

Indeed, let us assume the Lemma 5.1.2 true. Fixing $\xi$ and
$( \alpha_{a,b},f_{a,b})_{a,b}$ in (5.6), and letting
$(x, \omega )$ run in $\G$, gives:
\begin{equation}
|| \rho ( \xi ) ||_{\infty} \ \leq \
\sum_{b \in B} || \sum_{a  \in A_{b}} f_{a,b} ||_{\infty} .
\end{equation}
Then keeping $\xi$ fixed, but taking in (5.7) the infimum
after all the possible writings $\xi = \sum_{b \in B}
\sum_{a \in A_{b}} ( \alpha_{a,b},f_{a,b} )$, we obtain
\begin{equation}
|| \rho ( \xi ) ||_{\infty} \ \leq | \xi |_{o}^{'}.
\end{equation}
Finally, replacing in (5.8) $\xi $ by $\xi - \eta$, with
$\eta \in I( \Scal ),$ using that $\rho ( \eta ) =0$, and
taking another infimum, leads us to the inequality
\begin{equation}
|| \rho ( \xi ) ||_{\infty} \ \leq \ | \xi |_{o};
\end{equation}
this makes the inclusion $I_{o} ( \Scal ) \subseteq Ker \
\rho$ clear.

So we are left to make the

$\ $

{\bf Proof of Lemma 5.1.2} For every $b \in B$, denote
$\{ a \in A_{b} \ | \ \alpha_{a,b} \leq \beta_{x} \}$ by
$A_{b}^{'}$, and $A_{b} \setminus A_{b}^{'}$ by $A_{b}^{''}$.
By the definition of $\rho$ we have
$( \rho ( \xi ))(x, \omega )$ =
$\sum_{b \in B} \sum_{a \in A_{b}^{'}} f_{a,b} ( \omega )$,
which gives
$|( \rho ( \xi ))(x, \omega )| \ \leq$
$\sum_{b \in B} | \sum_{a \in A_{b}^{'}} f_{a,b} ( \omega )|$.
It remains to pick an element $b$ of $B$ and verify that
$| \sum_{a \in A_{b}^{'}} f_{a,b} ( \omega )| \ \leq$
$| \sum_{a \in A_{b}} f_{a,b} ( \omega )|$;
this is clear when $\sum_{a \in A_{b}^{'}} f_{a,b} =0$, so we
shall assume that there exists at least one $a \in A_{b}^{'}$
such that $f_{a,b}  ( \omega ) \neq 0$.
But, by an argument similar to the
one proving Lemma 5.1.1, it is seen that
$( \cup_{a \in A_{b}^{'}} \ supp \ f_{a,b} ) \ \cap$
$( \cup_{a \in A_{b}^{''}} \ supp \ f_{a,b} )$ =
$\emptyset$. Hence if there exists $a \in A_{b}^{'}$ such that
$f_{a,b} ( \omega ) \neq 0$, then we must have
$f_{a,b} ( \omega ) =0$ for all $a \in A_{b}^{''}$, so that
$| \sum_{a \in A_{b}^{'}} f_{a,b} ( \omega ) | $ =
$| \sum_{a \in A_{b}} f_{a,b} ( \omega ) |$.
${\bf QED}$

$\ $

{\bf 5.2 Remark} In Section 3.6 of \cite{K}, a class of
$\local$s with associated AF C*-algebras is constructed.
In particular, for a given sequence $(n_{i})_{i=0}^{\infty}$
of positive integers, the following $\local$
${\cal L}$ belonging to
this class has UHF($(n_{i})_{i=0}^{\infty}$) as associated
C*-algebra:

- the space on which ${\cal L}$ acts is
$X = \prod_{i=0}^{\infty} \{ 0,1, \ldots , n_{i}-1 \}$;

- ${\cal L}$ itself is described as
\begin{equation}
\{ \gamma ( u_{o}, \ldots , u_{k};v_{o}, \ldots ,v_{k} ) \ | \
k \geq 0, \ 0 \leq u_{j},v_{j} \leq n_{j}-1 \mbox{ for }
0 \leq j \leq k \},
\end{equation}
where for $u_{o}, \ldots ,u_{k},v_{o}, \ldots ,v_{k}$ as above
the domain of
$\gamma ( u_{o}, \ldots , u_{k};v_{o}, \ldots ,v_{k} )$ is
$\{ (w_{o},w_{1}$,$w_{2}$ , $\ldots ) \in X \ |$ $w_{j}=v_{j}$ for
$j \leq k \}$, its range is
$\{ (w_{o},w_{1},w_{2}, \ldots ) \in X \ |$
$w_{j}=u_{j}$ for $j \leq k \}$, and
$\gamma ( u_{o}, \ldots , u_{k};v_{o}, \ldots ,v_{k} )$
acts by replacing the first $k+1$ components of the sequence.

It is immediate that ${\cal L}$ of (5.10) is an $\f$-$\invsg$
(its maximal elements are those in (5.10)
having $u_{k} \neq v_{k}$); hence Theorem 5.1 applies. The
groupoid $\G$ associated as in 2.2 to the action of ${\cal L}$
on $X$ is easily seen to be (again) the Glimm groupoid. On
the other hand, it is obvious that ${\cal L}$ is not singly
generated, and that (although the acted space $X$ is the same
as in Remark 4.2) the intersection between ${\cal L}$ and the
$\star$-semigroup generated by the restricted odometer is
reduced to $\{ \epsilon , \theta \}$.

Thus actions of two rather different $\invsg$s
can give raise to the same groupoid (and hence the same
C*-algebra). In the present case, the reason which makes
this happen is that both the $\invsg$ constructions, and the
groupoid one, rely on the same method of finding m.a.s.a.'s
in AF-algebras (\cite{SV}, Section 1.1).

$\ $

{\bf 5.3 Remark} A $\local$ $\Scal$ on the space $\Omega$ is
called free (\cite{K}, Definition 7.1) if for every
$\alpha \in \Scal$ and $\omega \in Dom( \alpha )$ such that
$\alpha ( \omega ) = \omega$, there exists a neighborhood
$U$ of $\omega$ in $\Omega$ on which $\alpha$ acts as the
identity. If in addition $\Scal$ is assumed to be an
$\f$-inverse semigroup, it follows that every $\alpha \in \Scal$
which has a fixed point is an idempotent. Indeed,
$\alpha , \omega ,U$ being as above, we may assume
(by using (5.1)) that
the characteristic function $\chi$ of $U$ belongs to
$\Scal$.
{}From $\chi \leq \alpha$ we infer that
$\chi$ and $\alpha$ are majorized by the same maximal
element of $\Scal$, which can only be $\epsilon$, the unit;
but $\alpha \leq \epsilon$ means that $\alpha$ is idempotent.

As an immediate consequence, if $\Scal$ is a free $\local$ on
$\Omega$ and also an $\f$-inverse semigroup, then the groupoid $\G$
associated to $\Scal$ as in 2.2 is principal (\cite{Ren},
Definition I.1.1.1). It is actually clear that $\G$ coincides
with the equivalence relation ``$\omega \sim \omega '$
$\ecdef$ there is $\alpha \in \Scal$ such that
$\alpha ( \omega ) = \omega '$'' on $\Omega$; hence $\G$ is
exactly as in Section 7.2 of \cite{K}.

$\ $

{\bf 5.4 Example} An example of $\local$ which is an
$\f$-$\invsg$, but is not free, is the one described in
Section 10.1 of \cite{K}, which is related to the
Cuntz-Krieger algebras.

Recall from \cite{CK} that for an $n \times n$ matrix $A$
with entries $A_{i,j} \in \{ 0,1 \}$, $ 1 \leq i,j \leq n$,
which is irreducible (in the sense that for every $i,j$ there
exists $m \in \N$ such that $(A^{m})_{i,j} > 0$), and is not
a permutation matrix, there exists a unique C*-algebra
${\cal O}_{A}$ generated by $n$ non-zero partial isometries
$S_{1}, \ldots ,S_{n}$ satisfying
\begin{equation}
\left \{   \begin{array}{ll}
(S_{i}S_{i}^{*})(S_{j}S_{j}^{*}) \ =  \ 0,  &
\mbox{  for } i \neq j \\
S_{i}^{*} S_{i} \ =
\ \sum_{j=1}^{n} A_{i,j} (S_{j}S_{j}^{*}),  &
\mbox{  for } 1 \leq i \leq n .
\end{array}  \right.
\end{equation}

For the matrix $A$ as above, consider the compact space of
$A$-admissible sequences,
\[
X_{A} \ = \ \{ (j_{o},j_{1},j_{2}, \ldots ) \ | \
1 \leq j_{k} \leq n \mbox{ and } A_{j_{k},j_{k+1}} =1
\mbox{ for every } k \geq 0 \}
\]
(with topology obtained by restricting the product topology
of $\{ 0,1, \ldots ,n \}^{\N}$ ).
Then for every $1 \leq m \leq n$ denote by $\beta_{m}$ the map:
\[
\left\{  \begin{array}{l}
\{ (i_{k})_{k \geq 0} \in X_{A} \ | \ A_{m,i_{o}}=1 \}
\ \rightarrow  \
\{ (j_{k})_{k \geq 0} \in X_{A} \ | \ j_{o}=m \} \\
(i_{o},i_{1},i_{2}, \ldots )
\ \rightarrow  \ (m, i_{o},i_{1},i_{2}, \ldots ) ;
\end{array}  \right.
\]
and let $\Scal_{A}$ be the unital $\star$-semigroup of
$\homeo$s between open compact subsets of $X_{A}$ which is
generated by $\beta_{1}, \ldots , \beta_{n}$.
$\Scal_{A}$ is a $\local$; indeed, for every finite sequence
$j_{o}, \ldots ,j_{p}$ such that $A_{j_{o},j_{1}}= \ \ldots$
$=A_{j_{p-1},j_{p}}=1$, the domain of
$\beta_{j_{p}}^{*} \cdots \beta_{j_{o}}^{*}$ is the set
$V_{j_{o}, \ldots ,j_{p}}$ of sequences in $X_{A}$ beginning
with $j_{o}, \ldots ,j_{p}$, and the
$V_{j_{o}, \ldots ,j_{p}}$'s are a basis of $X_{A}$.
Note moreover that the domains of $\beta_{1}^{*}, \ldots ,
\beta_{n}^{*}$ form a partition of $X_{A}$, on which
$\beta_{1}^{*}, \ldots , \beta_{n}^{*}$ are the restrictions
of the one-sided shift on $X_{A}$; hence $\Scal_{A}$ is as
in 10.1 of \cite{K}.

On the other hand, it is not difficult to check that
$\Scal_{A}$ is an $\f$-$\invsg$;
we leave to the reader the verification of the following
facts:

$\ $

{\bf 5.4.1 Lemma} Let ${\bf {F_{n}}}$ be the free group on
generators $g_{1}, \ldots ,g_{n}$, and let $M_{A}$ be the
set of words $x= g_{i_{1}} \cdots g_{i_{p}} g_{j_{q}}^{-1}
\cdots g_{j_{1}}^{-1} \in {\bf {F_{n}}}$, with $p,q \geq 0$
and $1 \leq i_{1}, \ldots ,i_{p},j_{1}, \ldots ,j_{q} \leq n$,
such that:

(i) $i_{p} \neq j_{q}$ (i.e. the word $x$ is in reduced form);

(ii) $A_{i_{1},i_{2}}, \ldots ,A_{i_{p-1},i_{p}},
A_{j_{1},j_{2}}, \ldots ,A_{j_{q-1}.j_{q}}$ are all 1;

(iii) the set $\{ 1 \leq l \leq n \ | \
A_{i_{p},l} = A_{j_{q},l} =1 \}$ is not void.

For every $x= g_{i_{1}} \cdots g_{i_{p}} g_{j_{q}}^{-1}
\cdots g_{j_{1}}^{-1} \in M_{A}$ denote
$\beta_{i_{1}} \cdots \beta_{i_{p}} \beta_{j_{q}}^{*}
\cdots \beta_{j_{1}}^{*} \in \Scal_{A}$ by $\beta_{x}$.
Then the set of maximal elements of $\Scal_{A}$ is
$\{ \beta_{x} \ | \ x \in M_{A} \}$. Moreover, the
multiplicative structure on $M_{A}$ defined as in Section 1.5
above coincides with the one coming from ${\bf F_{n}}$ (i.e.
for $x,y \in M_{A}$ such that $\beta_{x} \beta_{y} \neq \theta$,
the product in ${\bf F_{n}}$ of $x$ and $y$ is still in
$M_{A}$, and $\beta_{x} \beta_{y} \leq \beta_{xy}$).

$\ $

{}From Theorem 5.1 above and from 10.1 of \cite{K} it follows
that the groupoid $\G_{A}$ associated (as in 2.2) to
$\Scal_{A}$ has C*($\G_{A}) \simeq {\cal O}_{A}$. Another
proof of this isomorphism can be done by using the
surjective $\homo$ $\Psi : \mbox{C*}( \Scal_{A} )
\rightarrow \mbox{C*}( \G_{A} )$ offered by Corollary 6.4
below: since $\beta_{1}, \ldots , \beta_{n}$ generate
$\Scal_{A}$, the partial isometries
$\Psi ( \beta_{1} ), \ldots , \Psi ( \beta_{n} )$
generate C*($\G_{A}$), and it is easily checked that the
latter ones satisfy (5.11) (thus C*$(\G_{A}) \simeq
{\cal O}_{A}$, due to the uniqueness of ${\cal O}_{A}$).

It should be noted that if the matrix $A$ has all the
entries equal to 1, then the groupoid $\G_{A}$ in the
preceding paragraph coincides with the Cuntz groupoid
of \cite{Ren}, Section III.2.1.

$\ $

\section{ The homomorphism C*($\Scal$)$\rightarrow$C*($\G$) }

Let $\Scal$ be an $\f$-inverse semigroup,
\setcounter{equation}{0}
and let $\Phi : \Scal \rightarrow {\cal I}_{\Omega}$ be an action
of $\Scal$ on a space $\Omega$, as in 2.1,
having the additional property that $Dom(\Phi (\alpha ))$ is compact
for every $\alpha \in \Scal$. (This is true for the examples
discussed in Sections 3, 4.2, 5.4.) Note that in particular
$\Omega = Dom(\Phi ( \epsilon ))$ must be compact.

Consider the groupoid $\G$ associated to the action $\Phi$
(as in 2.2). The partition
$\G = \cup_{x \in M} \{ x \} \times Dom(\Phi (\beta_{x}))$
remarked in 2.2 consists now of open and compact subsets;
each of these subsets is a $G$-set (
which means that the restriction to it of
the domain map and of the range map are one-to-one - see
\cite{Ren}, Definition I.1.10).

Recall now from \cite{Ren}, Definition I.2.10, that the open compact
$G$-sets of $\G$ form (with the pointwise multiplication) an inverse
semigroup, called the ample semigroup of $\G$. The partition mentioned
in the preceding paragraph (which is indexed by the set $\M$ of
maximal elements of $\Scal$) ``extends''
to a homomorphism of inverse semigroups as in the following Lemma, the
straightforward proof of which is left to the reader.

$\ $

{\bf 6.1 Lemma} For every $\alpha \in \Scal$ define
\begin{equation}
{\cal A} = \{ x \} \times Dom(\Phi (\alpha )) \subseteq \G ,
\end{equation}
where $x$ is the unique element of $M$ such that
$\alpha \leq \beta_{x}$. (If $\alpha$ has a zero element
$\theta$, we take by convention ${\cal A} (\theta) = \emptyset$,
the void set.) Then $\alpha \rightarrow {\cal A} (\alpha )$ is a unital
$\star$-homomorphism from $\Scal$ into the ample semigroup of $\G$.

$\ $

Now, it is a basic fact that if $A$ and $B$ are open compact $G$-sets
of an $r$-discrete
groupoid $\G$, and if $\chi_{A}$, $\chi_{B}$, $\chi_{AB}$
are the characteristic functions of $A, B, AB$, respectively, then
$\chi_{AB}$ is the convolution of $\chi_{A}$
and $\chi_{B}$ in $C_{c} (\G )$,
and moreover $\chi_{A}^{*} = \chi_{A^{-1}}$ in $C_{c} (\G )$.
As a consequence, Lemma 6.1 actually gives
a unital $\star$-homomorphism
$\Scal \rightarrow C_{c} (\G ) \subseteq C^{*} (\G )$.
This extends by linearity to the
``$\star$-semigroup algebra'' $\C [\Scal ]$.
(Recall that $\C [\Scal ]$ is the $\star$-algebra
having $\Scal$ (or $\Scal \setminus \{ \theta \}$, in
the case with zero element) as a linear basis, and with
multiplication and $\star$-operation coming
from those of $\Scal$.) Moreover, the unital $\star$-homomorphism
$\C [\Scal ] \rightarrow \mbox{C*} (\G )$
extends by universality to C*$(\Scal )$,
which is by definition the enveloping C*-algebra of $\C [\Scal ]$.

We pause here to recall that $\Scal$ has a left
regular representation (\cite{B}, p. 363);
its extension to $\C [\Scal ]$ is faithful (\cite{W}),
and this makes C*$(\Scal )$
to be a completion of $\C [\Scal ]$, and not of a
proper quotient of it. The envelopation of $\C [\Scal ]$
is done after all the unital $\star$-representations
of $\C [\Scal ]$
on Hilbert spaces, which correspond to $\star$-representations
by partial isometries of $\Scal$, and are all
automatically bounded with respect to the $l^{1}$-norm. (See
\cite{DP} for more details.)

Hence the Lemma 6.1 can be rephrased:

$\ $

{\bf 6.2 Proposition} There exists a unital  $\star$-homomorphism
$\Psi : \mbox{C*} (\Scal ) \rightarrow \mbox{C*} (\G )$,
uniquely determined by
\begin{equation}
\Psi (\alpha ) = \chi_{A(\alpha )},  \ \ \alpha \in \Scal
\end{equation}
(with $\A (\alpha )$ as in (6.1), and  $\chi_{A(\alpha )}
\in C_{c} (\G )$ its characteristic function).

$\ $

Let us next denote (as in the Example $1.4.1^{o}$ above) by
$\Scal^{(o)}$ the subsemigroup of idempotents of $\Scal$.
Then $clos \ sp \{ \gamma \ | \ \gamma \in \Scal^{(o)} \}$
$\subseteq \ \mbox{C*} (\Scal )$
is, clearly, a unital Abelian C*-subalgebra, which
will be denoted by C*$(\Scal^{(o)} )$.
(It is not difficult to show
that, actually, this really is canonically isomorphic
to the C*-algebra of the inverse semigroup $\Scal^{(o)}$.)
On the other hand we denote (following \cite{Ren},
Section II.4) by $C^{*} (\G^{(o)} )$
the unital Abelian C*-subalgebra
$\{ f \in C_{c} (\G ) \ | \ \ supp \ f \subseteq \G^{(o)}$ =
$\{ e \} \times \Omega \}$ of C*$(\G )$.
It is clear that $\Psi$ of (6.2) induces a unital $\star$-homomorphism
$\Psi^{(o)} : \mbox{C*} (\Scal^{(o)} )
\rightarrow \mbox{C*} (\G^{(o)} )$.

We have the following:

$\ $

{\bf 6.3 Theorem} In the context of Proposition 6.2,
assume in addition that $\Scal$ is
countable and that $\Omega$ (the space
of the action) is second countable. Then:

$1^{o}$ $\Psi$ is onto if and only if $\Psi^{(o)}$ is so.

$2^{o}$ $\Psi$ is an isomorphism if and only if $\Psi^{(0)}$ is so.

$\ $

{\bf Proof} If $\Psi$ is onto, then $\Psi^{(o)}$ is onto.
Indeed, it is known (\cite{Ren}, Proposition II.4.8, where we compose
on the right with the canonical surjection
$\mbox{C*} (\G ) \rightarrow {\mbox{C*}}_{red} (\G )$ )
that there exists a conditional expectation
$P: \mbox{C*} (\G ) \rightarrow \mbox{C*} (\G^{(o)} )$,
such that $P(f) = \chi_{\G^{(o)}}f$, $f \in C_{c} (\G )$.
For $\alpha \in \Scal \setminus \Scal^{(o)}$ it is clear that
$P(\Psi (\alpha )) = 0$, and this immediately implies
$Ran (P \circ \Psi ) \subseteq Ran(\Psi^{(o)})$;
but $P \circ \Psi$ is onto,
since $P$ and $\Psi$ are so, hence
$Ran(\Psi^{(o)}) = \mbox{C*} (\G^{(o)} )$.

It is obvious that $\Psi$ faithful $\Rightarrow  \ \Psi^{(o)}$
faithful, so only the parts
``$\Leftarrow $'' of $1^{o}$ and $2^{o}$ above remain to be discussed.
We shall use the following

$\ $

{\bf 6.3.1 Lemma} Let $x$ be in $M$  (the index set for the
maximal elements of $\Scal$) and consider the compact
open $G$-set $\A (\beta_{x}) = \{ x \} \times Dom(\Phi (\beta_{x}))$
of $\G$. Then $\Psi$ induces a contractive linear map
\begin{equation}
\Psi_{x} : clos \ sp \{ \alpha \in \Scal \ | \
\alpha \leq \beta_{x} \} \ \rightarrow \
\{ f \in C_{c} (\G ) \ | \ \ supp \ f \subseteq
\A (\beta_{x}) \} ,
\end{equation}
where the latter set is a closed linear subspace of C*$(\G )$.
If $\Psi^{(o)}$ is onto, then $\Psi_{x}$ is onto; if
$\Psi^{(o)}$ is an isomorphism, then $\Psi_{x}$ is an isomorphism.

$\ $

{\bf Proof of Lemma 6.3.1} If ${\cal C}$ is a C*-algebra ,
${\cal C}_{o} \subseteq {\cal C}$ is a C*-subalgebra,
and $w \in {\cal C}$
is a partial isometry such that $w^{*}w \in {\cal C}_{o}$,
then $w {\cal C}_{o}$ is a closed linear subspace of ${\cal C}$
(for instance because it can be written as
$\{ c \in {\cal C} \ | \ ww^{*} c = c, w^{*} c \in {\cal C}_{o} \}$).
Assume in addition that we also have
${\cal B}_{o} \subseteq {\cal B}$ C*-algebras,
$v \in {\cal B}$ partial isometry with $v^{*}v \in {\cal B}_{o}$,
and $\Psi : {\cal B} \rightarrow {\cal C}$ $\star$-homomorphism
such that $\Psi ( {\cal B}_{o} ) \subseteq {\cal C}_{o}$,
$\Psi (v)=w$. Then:
(a) $\Psi (v {\cal B}_{o} )  \subseteq w {\cal C}_{o}$;
(b) $\Psi ( {\cal B}_{o} ) = {\cal C}_{o}$ $\Rightarrow$
$\Psi (v {\cal B}_{o} ) =w {\cal C}_{o}$;
(c) $\Psi | {\cal B}_{o}$ faithful $\Rightarrow$
$\Psi | v{\cal B}_{o}$ faithful. Indeed, (a),(b) are clear, while
(c) comes out from the norm evaluation:
\[
|| \Psi (vb) ||^{2}  = || \Psi b^{*}(v^{*}v)b ||
\stackrel{( \star )}{=}
||b^{*}(v^{*}v)b||=||vb||^{2}, \ \ b \in {\cal B}_{o},
\]
(where at ($\star$) we used the faithfulness
of $\Psi | {\cal B}^{o}$).

In our context, we have to put ${\cal B}= \mbox{C*}(\Scal )$,
${\cal B}_{o} = \mbox{C*}(\Scal^{(o)} )$, $v= \beta_{x}$ and
${\cal C} = \mbox{C*}(\G )$,
${\cal C}_{o} = \mbox{C*}(\G^{(o)} )$,
$w= \chi_{\A ( \beta_{x} )}$ ($\Psi$ is $\Psi$).

$\ $

We return to the proof of the Theorem.
If $\Psi^{(o)}$ is onto, then from Lemma 6.3.1 and the direct
sum decomposition of (2.4) it is clear that
$Ran \Psi \supseteq C_{c} (\G )$, which is dense in C*$(\G )$;
hence $\Psi$ is onto.

If $\Psi^{(o)}$ is an isomorphism, then (by the Lemma)
every $\Psi_{x}$ of (6.3) is a linear isomorphism. The family of
linear maps $(\Psi^{-1}_{x})_{x \in M}$ defines (due to the
direct sum decomposition (2.4)) a linear map
$\rho : C_{c} (\G ) \rightarrow \mbox{C*} (\Scal )$. This is easily
checked to be a unital $\star$-homomorphism
of $C_{c} (\G )$ (one uses the corresponding properties of
$\Psi$). Now, due to the separability hypothesis made in the Theorem,
the groupoid $\G$ is second countable,  and we can
(exactly as we did in the proofs of the Theorems 4.1 and 5.1)
extend $\rho$ to a unital C*-algebra homomorphism
$\widetilde{\rho} : \mbox{C*} (\G ) \rightarrow \mbox{C*} (\Scal )$.
$\widetilde{\rho}$ clearly is an inverse for $\Psi$, which is
hence an isomorphism. {\bf QED}

$\ $

The homomorphism $\Psi^{(o)}$ of Theorem 6.3 has a natural
interpretation as a map between compact spaces. Indeed, let
$Z$ denote the spectrum of the (unital and Abelian)
C*-algebra C*$(\Scal^{(o)} )$; this is clearly
homomorphic (and will be henceforth identified) to the space of
multiplicative functions
$\zeta : \Scal^{(o)} \rightarrow \{ 0, 1\}$, such that
$\zeta (\epsilon ) = 1, \  \zeta (\theta ) = 0$
(where $\epsilon $ is the unit of $\Scal$, and $\theta $
its zero element - if it exists). On the other hand,
the convolution of the functions in C*$(\G^{(o)} ) $ coincides
with their pointwise product, which makes clear that the
spectrum of C*$(\G^{(o)} )$ is canonically identified to $\Omega $.

Every $\omega \in \Omega $ gives a character $\zeta_{\omega } \in Z$,
determined by $\zeta_{\omega } (\gamma ) = 1$, if
$\omega \in Dom( \Phi (\gamma ))$, and
$\zeta_{\omega} ( \gamma ) = 0$ otherwise
($\gamma \in \Scal^{(o)}$). The map
\begin{equation}
\psi^{(o)} : \Omega \rightarrow Z,  \ \
\psi^{(o)} (\omega ) = \zeta_{\omega }
\end{equation}
is continuous, because $Dom( \Phi (\gamma ))$
is open and compact for every
$\gamma \in \Scal^{(o)}$. It is an immediate verification (left to
the reader) that $\psi^{(o)}$ is the map between
character spaces corresponding to
$\Psi^{(o)} : \mbox{C*} (\Scal^{(o)} )
\rightarrow \mbox{C*} (\G^{(o)} )$
of Theorem 6.3. Hence we get:

$\ $

{\bf 6.4 Corollary} In the context of Theorem 6.3, the homomorphism
$\Psi : C^{*} (\Scal ) \rightarrow C^{*} (\G)$ is onto if and only
if $\psi^{(o)} : \Omega \rightarrow Z$ of (6.4)
is one-to one, i.e. if and only if  the subsets
$(Dom(\Phi (\alpha )))_{\alpha \in \Scal }$
separate the points of $\Omega $. Moreover, $\Psi$ is
an isomorphism if and only if $ \psi^{(o)}$ is bijective.

$\ $

{\bf 6.5 Example} Let us see how the Corollary 6.4 applies
in the context of Toeplitz inverse semigroups.
Consider $G$, $P$, $(\beta_{x})_{x \in G}$, $\sgp$ as
in Example 1.3, and let
$\Phi : \sgp \rightarrow {\cal I}_{\Omega }$ be the action defined
in Proposition 3.3. From equation (3.3) it is clear that
$Dom(\Phi (\alpha ))$
is open and compact in $\Omega $, for every $\alpha \in \sgp$,
hence the considerations made in this section
can be applied. Moreover, the family of subsets
$Dom(\Phi (\alpha ))_{\alpha \in \Scal_{G,P}}$
separates the points of $\Omega$. Indeed, for $A_{1} \neq A_{2}$
in $\Omega $ we can take an $x \in (A_{1} \setminus A_{2} ) \cup
(A_{2} \setminus A_{1} )
\subseteq \pq $ and it is obvious from equation (3.4) in Section 3
that $Dom(\Phi (\beta_{x}))$ will separate $A_{1}$ from $A_{2}$.

Hence, we can construct the natural homomorphism
$\Psi : \mbox{C*} (\sgp ) \rightarrow \mbox{C*} (\G )$,
determined by the equation (6.2), and moreover,
$\Psi$ is always surjective.

{}From Corollary 6.4 it also comes out that $\Psi$ will be an isomorphism
if and only if the natural map $\psi^{(0)}$ from $\Omega$ to
the space $Z$ of characters
of C*$(\Scal^{(o)}_{G,P} )$ is onto. Concerning this, we
mention without proof the following facts:

(a) $Z$ can be naturally identified to a subspace of $\{ 0,1 \}^{G}$,
in such a way that $\psi^{(o)} : \Omega \rightarrow Z$
becomes an inclusion. This comes from the fact that
$\Szero_{G,P}$ can be shown to be generated by the family
$(\beta_{x}^{*} \beta_{x} )_{x \in \pq}$;
hence a multiplicative function on $\Szero_{G,P}$ is determined by its
values on this family, and we get an embedding
$\tau : Z \rightarrow \{ 0,1 \}^{G}$,
defined by:
\begin{equation}
\tau (\zeta ) = \{ x \in \pq \ | \
\zeta ( \beta_{x} \beta_{x}^{*} ) = 1 \} .
\end{equation}
It turns out that $\Omega \stackrel{\psi^{(o)}}{\rightarrow} Z
\stackrel{\tau}{\rightarrow} \{ 0,1 \}^{G}$
is the inclusion of $\Omega$ into $\{ 0,1 \}^{G}$.

(b) It is well-known that if ${\cal E}$ is a unital semilattice,
then for every $\gamma \in {\cal E}$
which is not zero element we have a
character $\zeta_{\gamma }$ of C*$({\cal E})$, determined by:
$\zeta_{\gamma } (\gamma ') = 1$, if
$\gamma ' \geq \gamma$, $\zeta_{\gamma } (\gamma ') = 0$, if
$\gamma ' \not\geq \gamma$ $( \gamma '  \in {\cal E})$;
moreover, the family of the characters $( \zeta_{\gamma })_{\gamma }$
is dense in the spectrum of C*$({\cal E})$.
This offers the possibility of writing explicitly a dense subset of
$\tau (Z)$, i.e. of Z identified
inside $ \{ 0,1 \}^{G}$ as in the preceding
paragraph. The dense subset is
$\{ B_{x_{1},\ldots , x_{n}} \ |$
$n \geq 1, x_{1}, \ldots , x_{n} \in \pq \}$, where
\begin{equation}
B_{x_{1},\ldots , x_{n}}  \
\begin{array}[t]{l}
= \{ x \in \pq \ | \ \beta_{x} \beta_{x}^{*} \geq
(\beta_{x_{1}}^{*} \beta_{x_{1}})  \cdots
(\beta_{x_{n}}^{*} \beta_{x_{n}}) \}   \\
= \{ x \in \pq \ | \ xP \supseteq P \cap x^{-1}_{1}P \cap \dots
\cap  x^{-1}_{n}P \} .
\end{array}
\end{equation}

(c) Let us denote by ``$\prec$'' the left-invariant partial pre-order
determined by $P$ on $G$, i.e. $x \prec y \ecdef
x^{-1} y \in P$, for $x, y \in G$.
We shall call $(G, P)$ ``quasi-lattice ordered'' if
$P \cap P^{-1} = \{ e \}$
and if:

- every $x \in G$ having upper bounds in $P$ ( equivalently,
$x \in \pq$) has a least upper bound in $P$, denoted by $\sigma (x)$;

- every $s,t \in P$ having common upper bounds also have
a least common upper bound, denoted by  $\sigma (s,t)$.

The class of partially left-ordered groups satisfying these conditions
contains the totally left-ordered groups, and is closed under direct
products, semidirect products by order-preserving automorphisms, and
free products (see \cite{N2}), Example 2.3).

For $( G, P)$ quasi-lattice ordered and such that, say, $\pq \neq G$,
the Toeplitz inverse semigroup $\sgp$ is isomorphic to
$ (P \times P ) \cup \{ \theta \}$, with multiplication and
$*$-operation:
\begin{equation}
\left\{ \begin{array}{l}
(s,t) (u,v) = \left\{ \begin{array}{ll}
(s t^{-1} \sigma (t,u),v u^{-1} \sigma (t,u) ),  &
\mbox{if $t,u$ have common upper bounds} \\
\theta,   &  \mbox{otherwise}
\end{array} \right.  \\
\theta (s,t )  =  (s,t ) \theta = \theta = \theta^{*} = \theta^{2} \\
(s, t)^{*}  =  (t, s)
\end{array} \right.
\end{equation}
Indeed, $(s,t) \rightarrow \beta_{s} \beta_{t}^{*}$,
$\theta \rightarrow \theta$,
is easily seen to define an isomorphism between
$(P \times P ) \cup \{ \theta \}$
and $\sgp$. (In the case when $\pq = G$, we have a similar isomorphism
$P \times P \rightarrow \sgp$.)

If $(G, P)$ is quasi-lattice ordered then it is immediately verified,
using the particular form of
$\Scal_{G,P}$, that the sets in (6.6) belong
indeed to $\Omega$. Hence in this case we have a canonical isomorphism
C*$(\Scal_{G,P} ) \simeq  \mbox{C*} (\G )$.

(d) In general, the canonical $\homo$ of C*$(\Scal_{G,P} )$ onto
C*$(\G )$ is not faithful, as it can be seen on very simple examples
which don't have lattice properties. For instance for $G = \Z^{2}$
and $P = \{ (t_{1} t_{2} ) \in \Z^{2} \ |$
$0 \leq t_{2} \leq 2t_{1} \}$,
the identification of $Z$ inside $\{ 0,1 \}^{G}$ contains
elements which are not in $\Omega $, and which
are essentially produced by intersections of
lines $\{ t_{2} = a \} \cap \{ 2t_{1} - t_{2} = b \}$ with
$a, b \in \N$ of different parity.

$\ $

{\bf 6.6 Remark} Let $\Scal$ be an $\f$-$\invsg$. There
always exists a ``canonical'' action of
$\Scal$ which can be considered;
namely, $\Scal$ acts by conjugation on the spectrum of C*$(\Szero )$,
with $\Szero$ the subsemigroup of idempotents of
$\Scal$ (see \cite{P}, Section 3).
Indeed, $\Scal$ is first seen to act by conjugation on C*$(\Szero )$
(this being an action by not necessarily unital $\star$-endomorphisms);
passing to the space $Z$ of characters of C*$(\Szero  )$, we obtain an
action $\Phi : \Scal \rightarrow {\cal I}_{Z}$ described as follows:
\begin{equation}
\left\{ \begin{array}{ll}
Dom(\Phi (\alpha ))  =  \{ \zeta \in Z \ |
\ \zeta (\alpha^{*} \alpha ) = 1 \},  &  \ \alpha \in \Scal  \\
Ran(\Phi (\alpha ))  =  \{ \zeta \in Z \ |
\ \zeta (\alpha \alpha^{*} ) = 1 \},  &  \ \alpha \in \Scal \\
( \Phi (\alpha )) (\zeta )  (X) \ = \
\zeta ( \alpha^{*} X \alpha )   &  \
\alpha \in \Scal,  \zeta \in Dom( \Phi ( \alpha )),
X \in \mbox{C*} (\Szero ) .
\end{array} \right.
\end{equation}
It is clear that $Dom(\Phi (\alpha ))$
is an open and compact subset of $Z$,
for every $\alpha \in \Scal$, hence the
considerations made in this section can
be applied. Moreover, the map $\psi^{(0)}$
discussed in Corollary 6.4 clearly
is, in this case, the identity map of $Z$.
Hence if $\Scal$ is countable,
Corollary 6.4 gives that the associated groupoid $\G$ has
C*$(\G ) \simeq  \mbox{C*} (\Scal )$.

\pagebreak


\begin{thebibliography}{99}

\bibitem{B} B.A. Barnes.
Representations of the $l^{1}$-algebra of an inverse
semigroup, Trans. Amer. Math. Soc. 218(1976), 361-396.

\bibitem{CK} J. Cuntz, W. Krieger.
A class of C*-algebras and topological Markov chains,
Inventiones Math. 56(1980), 251-268.

\bibitem{DP} J. Duncan, A.L.T. Paterson.
C*-algebras of $\invsg$s, Proc. of the Edinburgh Math.
Soc. 28(1985), 41-58.

\bibitem{E1} R. Exel.
Circle actions on C*-algebras, partial automorphisms and
a generalized Pimsner-Voiculescu exact sequence, preprint.

\bibitem{E2} R. Exel.
Approximately finite C*-algebras and partial automorphisms,
preprint.

\bibitem{K} A. Kumjian.
On localizations and simple C*-algebras, Pacific J. Math.
112(1984), 141-192.

\bibitem{MR} P. Muhly, J. Renault.
C*-algebras of multivariable Wiener-Hopf operators, Trans.
Amer. Math. Soc. 274(1982), 1-44.

\bibitem{N1} A. Nica.
Some remarks on the groupoid approach to Wiener-Hopf operators,
J. Operator Theory 18(1987), 163-198.

\bibitem{N2} A. Nica.
C*-algebras generated by isometries and Wiener-Hopf operators,
Preprint INCREST No.33/1989, to appear in the Journal of Operator
Theory.

\bibitem{P} A.L.T. Paterson.
Inverse semigroups, groupoids, and a problem of J. Renault,
preprint.

\bibitem{Pet} M. Petrich.
Inverse semigroups, John Wiley, 1984.

\bibitem{Ren} J. Renault.
A groupoid approach to C*-algebras, Lecture Notes in Math. 793,
Springer Verlag 1980.

\bibitem{R} J. Renault.
Representation des produits croises
d'algebres de groupoides, J. Operator Theory 18(1987),
67-97.

\bibitem{SV} S. Stratila, D. Voiculescu.
Representations of AF-algebras and of the group $U( \infty )$,
Lecture Notes in Math. 486, Springer Verlag 1975.

\bibitem{W} J.R. Wordingham.
The left regular $\star$-representation of an $\invsg$,
Proc. Amer. Math. Soc. 86(1982), 55-58.

\end{thebibliography}
\end{document}